\documentclass[a4paper,11pt,leqno]{article}
\usepackage{amsmath,amssymb,latexsym,amscd,amsthm,amsxtra,graphicx}
\begin{document}




\newtheorem{thm}{Theorem}[subsection]
\newtheorem{lem}[thm]{Lemma}
\newtheorem{cor}[thm]{Corollary}
\newtheorem{prop}[thm]{Proposition}
\newtheorem{clm}{Claim}
\newtheorem{dfn}[thm]{Definition}


\newcommand{\pf}{{\it Proof} \,\,\,\,\,\,}
\newcommand{\remark}{{\it Remark} \,\,\,\,}

\newcommand{\cpx}{{\mathbb C}}
\newcommand{\rel}{{\mathbb R}}
\newcommand{\rat}{{\mathbb Q}}
\newcommand{\itg}{{\mathbb Z}}
\newcommand{\nat}{{\mathbb N}}
\newcommand{\qrt}{{\mathbb H}}
\newcommand{\qtr}{{\mathbb H}}
\newcommand{\dbP}{{\mathbb P}}
\newcommand{\dbS}{{\mathbb S}}
\newcommand{\dbO}{{\mathbb O}}
\newcommand{\dbK}{{\mathbb K}}

\newcommand{\bfa}{\mbox{\underline{\bf a}}}
\newcommand{\bfb}{\mbox{\underline{\bf b}}}
\newcommand{\bfc}{\mbox{\underline{\bf c}}}
\newcommand{\bfd}{\mbox{\underline{\bf d}}}
\newcommand{\bfe}{\mbox{\underline{\bf e}}}
\newcommand{\bff}{\mbox{\underline{\bf f}}}
\newcommand{\bfg}{\mbox{\underline{\bf g}}}
\newcommand{\bfh}{\mbox{\underline{\bf h}}}
\newcommand{\bfi}{\mbox{\underline{\bf i}}}
\newcommand{\bfj}{\mbox{\underline{\bf j}}}
\newcommand{\bfk}{\mbox{\underline{\bf k}}}
\newcommand{\bfl}{\mbox{\underline{\bf l}}}
\newcommand{\bfm}{\mbox{\underline{\bf m}}}
\newcommand{\bfn}{\mbox{\underline{\bf n}}}
\newcommand{\bfo}{\mbox{\underline{\bf o}}}
\newcommand{\bfp}{\mbox{\underline{\bf p}}}
\newcommand{\bfq}{\mbox{\underline{\bf q}}}
\newcommand{\bfr}{\mbox{\underline{\bf r}}}
\newcommand{\bfs}{\mbox{\underline{\bf s}}}
\newcommand{\bft}{\mbox{\underline{\bf t}}}
\newcommand{\bfu}{\mbox{\underline{\bf u}}}
\newcommand{\bfv}{\mbox{\underline{\bf v}}}
\newcommand{\bfw}{\mbox{\underline{\bf w}}}
\newcommand{\bfx}{\mbox{\underline{\bf x}}}
\newcommand{\bfy}{\mbox{\underline{\bf y}}}
\newcommand{\bfz}{\mbox{\underline{\bf z}}}
\newcommand{\bfmu}{{\underline{\bf \mu}}}

\newcommand{\bpa}{\mbox{\underline{\bf a}}}
\newcommand{\bpb}{\mbox{\underline{\bf b}}}
\newcommand{\bpc}{\mbox{\underline{\bf c}}}
\newcommand{\bpd}{\mbox{\underline{\bf d}}}
\newcommand{\bpe}{\mbox{\underline{\bf e}}}
\newcommand{\bpf}{\mbox{\underline{\bf f}}}
\newcommand{\bpg}{\mbox{\underline{\bf g}}}
\newcommand{\bph}{\mbox{\underline{\bf h}}}
\newcommand{\bpi}{\mbox{\underline{\bf i}}}
\newcommand{\bpj}{\mbox{\underline{\bf j}}}
\newcommand{\bpk}{\mbox{\underline{\bf k}}}
\newcommand{\bpl}{\mbox{\underline{\bf l}}}
\newcommand{\bpm}{\mbox{\underline{\bf m}}}
\newcommand{\bpn}{\mbox{\underline{\bf n}}}
\newcommand{\bpo}{\mbox{\underline{\bf o}}}
\newcommand{\bbpp}{\mbox{\underline{\bf p}}}
\newcommand{\bpq}{\mbox{\underline{\bf q}}}
\newcommand{\bpr}{{\mbox{\underline{\bf r}}}}
\newcommand{\bps}{\mbox{\underline{\bf s}}}
\newcommand{\bpt}{\mbox{\underline{\bf t}}}
\newcommand{\bpu}{\mbox{\underline{\bf u}}}
\newcommand{\bpv}{\mbox{\underline{\bf v}}}
\newcommand{\bpw}{\mbox{\underline{\bf w}}}
\newcommand{\bpx}{\mbox{\underline{\bf x}}}
\newcommand{\bpy}{\mbox{\underline{\bf y}}}
\newcommand{\bpz}{\mbox{\underline{\bf z}}}

\newcommand{\eps}{\varepsilon}

\newcommand{\zzz}{{\mathfrak z}}
\newcommand{\fff}{{\mathfrak f}}
\newcommand{\www}{{\mathfrak w}}
\newcommand{\vvv}{{\mathfrak v}}
\newcommand{\gggg}{{\mathfrak g}}
\newcommand{\kkk}{{\mathfrak k}}
\newcommand{\aaa}{{\mathfrak a}}
\newcommand{\saa}{\mbox{}^s{\mathfrak a}}
\newcommand{\aas}{{\mathfrak a}_S}
\newcommand{\bbb}{{\mathfrak b}}
\newcommand{\sbb}{\mbox{}^s{\mathfrak p}}
\newcommand{\spp}{\mbox{}^s{\mathfrak p}}
\newcommand{\llll}{{\mathfrak l}}
\newcommand{\tll}{\tilde{\mathfrak l}}
\newcommand{\eee}{{\mathfrak e}}
\newcommand{\ccc}{{\mathfrak c}}
\newcommand{\ttt}{{\mathfrak t}}
\newcommand{\stt}{\mbox{}^s{\mathfrak t}}
\newcommand{\lls}{{\mathfrak l}_S}
\newcommand{\obb}{\bar{\mathfrak b}}
\newcommand{\aad}{{\mathfrak a}^\ast}
\newcommand{\ads}{{\mathfrak a}^\ast_S}
\newcommand{\nnn}{{\mathfrak n}}
\newcommand{\nns}{{\mathfrak n}_S}
\newcommand{\uuu}{{\mathfrak u}}
\newcommand{\tuu}{\tilde{\mathfrak u}}
\newcommand{\buu}{\bar{\mathfrak u}}
\newcommand{\bnn}{\bar{\mathfrak n}}
\newcommand{\bab}{\bar{\mathfrak b}}
\newcommand{\bns}{\bar{\mathfrak n}_S}
\newcommand{\mmm}{{\mathfrak m}}
\newcommand{\smm}{\mbox{}^s{\mathfrak m}}
\newcommand{\mms}{{\mathfrak m}_S}
\newcommand{\jjj}{{\mathfrak j}}
\newcommand{\hhh}{{\mathfrak h}}
\newcommand{\shh}{\mbox{}^s{\mathfrak h}}
\newcommand{\uhh}{\mbox{}^u{\mathfrak h}}
\newcommand{\hhd}{{\mathfrak h}^\ast}
\newcommand{\qqq}{{\mathfrak q}}
\newcommand{\ppp}{{\mathfrak p}}
\newcommand{\tpp}{\tilde{\mathfrak p}}
\newcommand{\tnn}{\tilde{\mathfrak n}}
\newcommand{\tqq}{\tilde{\mathfrak q}}
\newcommand{\pps}{{\mathfrak p}_S}
\newcommand{\bapp}{\bar{\mathfrak p}}
\newcommand{\baps}{\bar{\mathfrak p}_S}
\newcommand{\sss}{{\mathfrak s}}
\newcommand{\ooo}{{\mathfrak o}}
\newcommand{\ddd}{{\mathfrak d}}

\newcommand{\ggc}{{\mathfrak g}_{\sf c}}
\newcommand{\bbc}{{\mathfrak b}_{\sf c}}
\newcommand{\llc}{{\mathfrak l}_{\sf c}}
\newcommand{\uuc}{{\mathfrak u}_{\sf c}}
\newcommand{\nsc}{{({\mathfrak n}_S)}_{\sf c}}
\newcommand{\nnc}{{({\mathfrak n}_S)}_{\sf c}}
\newcommand{\buc}{\bar{\mathfrak u}_{\sf c}}
\newcommand{\bnc}{{(\bar{\mathfrak n}_S)}_{\sf c}}
\newcommand{\hhc}{{\mathfrak h}_{\sf c}}
\newcommand{\hcd}{{{\mathfrak h}_{\sf c}}^\ast}
\newcommand{\kkc}{{\mathfrak k}_{\sf c}}
\newcommand{\ppc}{{\mathfrak p}_{\sf c}}
\newcommand{\pcs}{({\mathfrak p}_S)_{\sf c}}
\newcommand{\ncs}{({\mathfrak n}_S)_{\sf c}}
\newcommand{\bpp}{\bar{\mathfrak p}}
\newcommand{\bqq}{\bar{\mathfrak q}}

\newcommand{\spl}{\mbox{}^s}

\newcommand{\gl}{{\mathfrak g}{\mathfrak l}}
\newcommand{\gll}{{\mathfrak gl}_L}
\newcommand{\glr}{{\mathfrak gl}_R}
\newcommand{\so}{\sss\ooo}

\newcommand{\sH}{{\mbox{}^s{H}}}

\newcommand{\hol}{{\cal O}}
\newcommand{\dif}{{\cal D}}
\newcommand{\ana}{{\cal A}}

\newcommand{\tpca}{\tilde{\cal P}}
\newcommand{\aca}{{\cal A}}
\newcommand{\bca}{{\cal B}}
\newcommand{\cca}{{\cal C}}
\newcommand{\dca}{{\cal D}}
\newcommand{\eca}{{\cal E}}
\newcommand{\fca}{{\cal F}}
\newcommand{\gca}{{\cal G}}
\newcommand{\hca}{{\cal H}}
\newcommand{\ica}{{\cal I}}
\newcommand{\jca}{{\cal J}}
\newcommand{\kca}{{\cal K}}
\newcommand{\lca}{{\cal L}}
\newcommand{\mca}{{\cal M}}
\newcommand{\nca}{{\cal N}}
\newcommand{\oca}{{\cal O}}
\newcommand{\pca}{{\cal P}}
\newcommand{\qca}{{\cal Q}}
\newcommand{\rca}{{\cal R}}
\newcommand{\sca}{{\cal S}}
\newcommand{\tca}{{\cal T}}
\newcommand{\uca}{{\cal U}}
\newcommand{\vca}{{\cal V}}
\newcommand{\wca}{{\cal W}}
\newcommand{\xca}{{\cal X}}
\newcommand{\yca}{{\cal Y}}
\newcommand{\zca}{{\cal Z}}

\newcommand{\gS}{{\mathfrak S}}
\newcommand{\ii}{\sqrt{-1}}
\newcommand{\real}{\mbox{Re}}
\newcommand{\Real}{\mbox{Re}}
\newcommand{\res}{\mbox{res}}
\newcommand{\supp}{\mbox{supp}}
\newcommand{\rank}{\mbox{rank}}
\newcommand{\card}{\mbox{card}}
\newcommand{\Ad}{\mbox{Ad}}
\newcommand{\ad}{\mbox{ad}}
\newcommand{\dg}{{\deg}}
\newcommand{\Hom}{\mbox{Hom}}
\newcommand{\End}{\mbox{End}}
\newcommand{\Dim}{\mbox{Dim}}
\newcommand{\gr}{\mbox{gr}}
\newcommand{\Ass}{\mbox{Ass}}
\newcommand{\WF}{\mbox{WF}}
\newcommand{\Ann}{\mbox{Ann}}
\newcommand{\tr}{\mbox{tr}}
\newcommand{\Mod}{{\sf Mod}}
\newcommand{\sgn}{\mbox{sgn}}
\newcommand{\Ind}{\mbox{\sf Ind}}
\newcommand{\SO}{\mbox{SO}}
\newcommand{\Oo}{\mbox{O}}
\newcommand{\SOo}{{\mbox{SO}_0}}
\newcommand{\GL}{\mbox{GL}}
\newcommand{\Spp}{\mbox{Sp}}
\newcommand{\U}{\mbox{U}}
\newcommand{\triv}{{\mbox{triv}}}
\newcommand{\Rea}{{\mbox{Re}}}
\newcommand{\diag}{\mbox{diag}}


\title{{\sf The edge-of-wedge type embeddings of derived functor modules
for the type A classical groups}}
\author{{\sf Hisayosi Matumoto}\thanks{Supported in part by Grant-in-aid for Scientific Research (No.\ 50272597) .} 
\\Graduate School of Mathematical Sciences\\ University of Tokyo\\ 3-8-1
Komaba, Tokyo\\ 153-8914, JAPAN\\ e-mail: hisayosi@ms.u-tokyo.ac.jp}
\date{}
\maketitle
\begin{abstract}
A holomorphic discrete series  can be realized as the space of the
 holomorphic sections of a homogeneous vector bundle on a bounded
 symmetric domain.
We can embed  it into a degenerate principal series realized as the
 space of hyperfunction sections of a vector bundle on the Shilov
 boundary by taking the boundary value.
We consider an algebraic version of higher cohomological analogue of such
 an embedding for complex reductive groups and classical groups of type A.
\end{abstract}
\setcounter{section}{0}
\setcounter{subsection}{0}

\setcounter{section}{1}
\section*{\S\,\, 0.\,\,\,\, Introduction}

Let $G$ be a real linear reductive Lie group and let $G_\cpx$ its
complexification.
We denote by $\gggg_0$ (resp.\ $\gggg$) the Lie algebra of $G$ (resp.\
$G_\cpx$) and denote by $\sigma$ the complex conjugation on $\gggg$ with
respect to $\gggg_0$.
We fix a maximal compact subgroup $K$ of $G$ and denote by $\theta$ the
corresponding Cartan involution. 
We denote by $\kkk$ the complexified Lie algebra of $K$.

We fix a parabolic subgroup $P$ of $G$ with a $\theta$-stable Levi part
$M$.
We denote by $N$ the nilradical of $P$.
We denote by $\ppp$, $\mmm$, and $\nnn$  the complexified Lie algebras
of $P$, $M$, and $N$, respectively.
We denote by $P_\cpx$, $M_\cpx$, and $N_\cpx$ the analytic subgroups in
$G_\cpx$ corresponding to $\ppp$, $\mmm$, and $\nnn$, respectively. 

For $X\in \mmm$, we define 
\[
 \delta(X)=\frac{1}{2}\tr\left(\ad_\gggg(X)|_\nnn\right).
\]
Then, $\delta$ is a one-dimensional representation of $\mmm$.
We see that $2\delta$ lifts to a holomorphic group homomorphism $\xi_{2\delta}
: M_\cpx\rightarrow\cpx^\times$.
Defining $\xi_{2\delta}|_{N_\cpx}$ trivial, we may extend $\xi_{2\delta}
 $ to $P_\cpx$.
We put $X=G_\cpx/P_\cpx$.
Let $\lca$ be the holomorphic line bundle on $X$ corresponding to the
canonical divisor.
Namely, $\lca$ is the $G_\cpx$-homogeneous line bundle on $X$ associated
to the character $\xi_{2\delta}$ on $P_\cpx$. We denote the restriction
of $\xi_{2\delta}$  to $P$ by the same letter.

For a character $\eta : P\rightarrow \cpx^\times$, we consider the
unnormalized parabolic induction ${}^u\Ind_P^G(\eta)$.
Namely, ${}^u\Ind_P^G(\eta)$ is the $K$-finite part of the space of the
$C^\infty$-sections of the $G$-homogeneous line bundle on $G/P$
associated to $\eta$. ${}^u\Ind_P^G(\eta)$ is a Harish-Chandra $(\gggg, K)$-module.

If $G/P$ is orientable, then the trivial $G$-representation is the
unique irreducible quotient of ${}^u\Ind_P^G(\xi_{2\delta})$.
If $G/P$ is not orientable, there is a character $\omega$ on $P$ such
that $\omega$ is trivial on the identical component of $P$ and the trivial $G$-representation is the
unique irreducible quotient of
${}^u\Ind_P^G(\xi_{2\delta}\otimes\omega)$.

Let $\hol$ be an open $G$-orbit on $X$.
We put the following assumption:

{\bf Assumption A} \,\,\,\,\, There is a $\theta$-stable parabolic
subalgebra $\qqq$ of $\gggg$ such that $\qqq\in\hol$.

Under the above assumption, $\qqq$ has a Levi decomposition
$\qqq=\llll+\uuu$ such that $\llll$ is a $\theta$ and $\sigma$-stable
Levi part. 
In fact $\llll$ is unique, since we have $\llll=\sigma(\qqq)\cap\qqq$. 

For each open $G$-orbit $\hol$ on $X$, we put
\[
 \ana_\hol= \mbox{H}^{\dim \uuu\cap\kkk}(\hol, \lca)_{\mbox{$K$-finite}}.
\]
From [Wong 1993], $\ana_\hol$ is a derived functor module.
Namely, in the terminology in [Vogan-Zuckerman 1984], we have
$\ana_\hol=\ana_\qqq=\ana_\qqq(0)$.

If $G$ is of the Hermitian type and $P$ is a Siegel parabolic
subgroup of $G$,
there are two open $G$-orbits in $X$ isomorphic to $G/K$. (the Siegel
upper and lower half planes.)
For such orbits, $\ana_\hol$ are representations in holomorphic or anti-holomorphic discrete
series.
In this case, we can embed $\ana_\hol$  into ${}^u\Ind_P^G(\xi_{2\delta})$ by taking boundary
values at the Shilov boundary. 

We consider the other orbits.
Then, the corresponding embedding should be a higher
cohomological analogy of boundary value maps.
Here, we discuss examples.

{\it Example 1.} \,\,\,\,
Let $n$ be an integer such that $n\geqslant 3$.
Let $G$ be $SO_0(n,1)$ and let $P$ be its minimal parabolic subgroup.
Then $X=G_\cpx/P_\cpx$  has two $G$-orbits.
One is a closed orbit $X_\rel=G/P$ and the other is an open orbit (say
$\hol$).
In this case,  Assumption A holds and we have $\hol\cong G/(SO(2)\times SO_0(n-2,1))$.
Since $\hol=X-X_\rel$, we have the following long exact sequence of
cohomologies:
\[
 H^{n-2}(X,\lca)\rightarrow H^{n-2}(\hol,\lca)\rightarrow
 H^{n-1}_{X_\rel}(X,\lca)
\rightarrow H^{n-1}(X,\lca)\rightarrow H^{n-1}(\hol,\lca).
\]
Taking account of $n-1=\dim X$, we see $H^{n-2}(X,\lca)=0$ and
$H^{n-1}(X,\lca)=\cpx$ from the
generalized Borel-Weil-Bott theorem [Kostant 1961].
From [Wong 1993], we have
$H^{n-2}(\hol,\lca)_{\mbox{$K$-finite}}=\ana_\hol$ and $H^{n-1}(\hol,\lca)=0$.
Since $X_\rel$ is an $n-1$-dimensional sphere, $X_\rel$ is orientable.
Hence $H^{n-1}_{X_\rel}(X,\lca)$ is nothing but the space of
hyperfunction sections of the principal series representation of $G$ with respect to
$\lca$. (cf.\ [Kashiwara-Kawai-Kimura 1986])
Hence, taking $K$-finite parts, we have the following short exact
sequence.
\[
 0\rightarrow\ana_\hol\rightarrow {}^uInd_P^G(\xi_{2\delta})\rightarrow
 \cpx\rightarrow 0.
\]
Hence, the structure of ${}^uInd_P^G(\xi_{2\delta})$ is understood geometrically
in this case.

{\it Example 2.}  \,\,\,\,
Let $n\geqslant 3$.
Let $G=SO_0(n,2)$ and $P$ be a parabolic subgroup with the semisimple
part of the  Levi part is
isomorphic to $SO_0(n-1,1)$.
In this case $X=G_\cpx/P_\cpx$ has three open $G$-orbits.
Two of them (say $\hol_+$ and $\hol_-$) are Hermitian symmetric spaces
(symmetric domains of type IV).
The remaining one (say $\hol_0$)is non-Stein and isomorphic to
$G/(SO(2)\times SO_0(n-2,2))$.
Let $\overline{\hol}_+$  and $\overline{\hol}_-$ be the closures of
$\hol_+$ and $\hol_-$, respectively.
In this case, we have
$(X-\overline{\hol}_+)\cap(X-\overline{\hol}_-)=\hol_0$ and
$(X-\overline{\hol}_+)\cup(X-\overline{\hol}_-)=X-X_\rel$.
Here, we put $X_\rel=G/P$.
Hence, we have the following Mayer-Vietris exact sequence:
\[ 
H^{2n-2}(X-\overline{\hol}_+,\lca)\oplus H^{2n-2}(X-\overline{\hol}_-,\lca)\rightarrow H^{2n-2}(\hol_0,\lca)\rightarrow H^{2n-1}(X-X_\rel,\lca).
\]
We also have the following exact sequences.
\begin{align*} 
H^{2n-2}(X, \lca)\rightarrow H^{2n-2}(X-\overline{\hol}_+,\lca)\rightarrow H^{2n-1}_{\overline{\hol}_+}(X, \lca),\\
H^{2n-2}(X, \lca)\rightarrow
 H^{2n-2}(X-\overline{\hol}_-,\lca)\rightarrow
 H^{2n-1}_{\overline{\hol}_-}(X, \lca).
\end{align*}
From [Kostant 1961], we have $H^{2n-2}(X, \lca)=0$.
We can regard $\overline{\hol}_+$ and $\overline{\hol}_-$ as
closed convex sets in an open cell of $X$.
(Namely,  we  realize  $\hol_+$ and $\hol_-$ as bounded symmetric domains.)
Hence, from the edge-of-wedge theorem
(see [Kashiwara-Laurent 1983] Th\'{e}or\`{e}me 1.1.2),
we have $H^{2n-1}_{\overline{\hol}_+}(X,
\lca)=H^{2n-1}_{\overline{\hol}_-}(X, \lca)=0$.
Hence, $H^{2n-2}(X-\overline{\hol}_+,\lca)=H^{2n-2}(X-\overline{\hol}_-,\lca)=0.$
We also have the following exact sequence. 
\[
0=H^{2n-1}(X, \lca)\rightarrow  H^{2n-1}(X-X_\rel,\lca) \rightarrow H^{2n}_{X_\rel}(X,\lca).
\]
Hence, we have
\[
 H^{2n-2}(\hol_0,\lca)\hookrightarrow
 H^{2n-1}(X-X_\rel,\lca)\hookrightarrow H^{2n}_{X_\rel}(X,\lca).
\]
In this case $X_\rel$ is orientable if and only if $n$ is even.
Hence, the local cohomology $H^{2n}_{X_\rel}(x,\lca)$ is the space of
hyperfunction sections of the degenerate principal series of $G$ with
respect to $\lca$ (resp.\ $\lca\otimes\omega$) if $n$ is even (resp.\
odd).
Taking $K$-finite part, we have
\begin{align*}
\ana_{\hol_0}\hookrightarrow {}^u\Ind^G_P(\xi_{2\delta})  \,\,\,\,
 (\mbox{if $n$ is even}),\\
\ana_{\hol_0}\hookrightarrow {}^u\Ind^G_P(\xi_{2\delta}\otimes\omega)  \,\,\,\,
 (\mbox{if $n$ is odd}). 
\end{align*}

Such an easy construction of an embedding seems to be difficult to obtain in
the general cases.
However, [Matumoto 1988],[Kashiwara-Schapira 1990], [Sahi-Stein 1992], [Gindikin 1993],
etc. give some evidences of the existence of the edge-of-wedge type
embeddings as above
in more general settings.
In this article, we consider the following problem in the setting of
Harish-Chandra modules:

{\bf Problem B} \,\,\,\, Let $G$ be a real linear reductive Lie group
and let $P$ be its parabolic subgroup. 
Let $\hol$ be an open $G$-orbit in $X=G_\cpx/P_\cpx$ satisfying
Assumption A.
Does there exist a character $\chi$ of $P$ trivial on the
 identity component such that 
  $\ana_\hol \hookrightarrow {}^u\Ind_P^G(\xi_{2\delta}\otimes\chi)$?

If the nilradical of $P$ is commutative, the answer is known more or
less by [Sahi
1992], [Sahi-Stein 1990], [Sahi 1995], and [Zhang 1995].
For instance, we consider a real rank $n$ group $G=Sp(n,\rel)$ and its
Siegel parabolic subgroup $P$.
If $n$ is even, then Problem B is affirmative.
However, if $n$ is odd, Problem B fails except for holomorphic and
anti-holomorphic discrete series.

In this article, we first consider the case that $G$ is a complex
reductive group.
In this case, $X=G_\cpx/P_\cpx=G/P\times G/P$ has a unique open
$G$-orbit (say $\hol_0$).
Let $\gggg=Lie(G)$ and $\ppp=Lie(P)$.
We denote by $\mmm$ a Levi part of $\ppp$.
We fix a Cartan subalgebra $\hhh$
of $\gggg$ such that $\hhh\subseteq\ppp$.
We denote by $w_0$ (resp.\ $w_\ppp$) the longest element of the Weyl
group with respect to $(\gggg,\hhh)$  (resp.\ $(\mmm, \hhh)$).
$\hol_0$ satisfies the Assumption A if and only if
 $w_0 w_\ppp$ is an involution.

In fact, we have :

{\bf Theorem C} \,\,\, (Theorem 1.2.6)
{\it Assume $G$ is a complex reductive group.
$\ana_{\hol_0}\hookrightarrow {}^u\Ind^{G}_{P}(\xi_{2\delta} )$ if and only if  $w_0 w_\ppp$ is a Duflo involution in the Weyl group for
 $(\gggg,\hhh)$.
}

The condition `` $w_0 w_\ppp$ is a Duflo involution.'' is same
as the existence condition for a non-trivial homomorphism between
certain generalized Verma modules ([Matumoto 1993]).
In fact, $\ana_{\hol_0}$ can be identified with some non-class one
degenerate principal series representation ([Enright 1979] and [Vogan-Zuckerman 1984]).
In viewpoint of this identification, the embedding in Theorem C is
nothing but an intertwining operator between degenerate principal seiries
induced from the homomorphism between generalized Verma modules. 
(cf.\ [Collingwood-Shelton 1990])

For a Weyl group of type A, any involution is a Duflo involution.
We get an affirmative answer on the case of $G=\GL(n,\cpx)$. 
Moreover, we have:

{\bf Theorem D} \,\,\, (Theorem 3.2.1)
{\it Under Assumption A, we have an affirmative answer on the case of
$G=\GL(n,\rel)$ and $G=\GL(n,\qtr)$.
}

In this case, as a submodule of
${}^u\Ind^{G_\cpx}_{P_\cpx}(\xi_{2\delta} )$, $\ana_{\hol_0}$ is the image of intertwining
operator from the degenerate principal series representation of which $\ana_{\hol_0}$
is realized as a quotient.
This intertwing operator is induced from a homomorphism between
generalized Verma modules. 

Finally, we mention the case of $G=\U(m,n)$.

Let $G=\U(m,n)$ and let $P$ be an arbitrary prarabolic subgroup of $G$.
In this case, Assumption A automatically holds for all the open
orbits.
However, there is an obvious necessary condition.
Since ${}^u\Ind^{G}_{P}(\xi_{2\delta} )$ has a perfect pairing
with an irreducible generalized Verma module, any  submodule of
${}^u\Ind^{G}_{P}(\xi_{2\delta} )$ has the same
Gelfand-Kirillov dimension as ${}^u\Ind^{G}_{P}(\xi_{2\delta}
)$.
We call an open $G$-orbit $\hol$ in $X=G_\cpx/P_\cpx$ good, if the
Gelfand-Kirillov dimension of $\ana_{\hol}$ equals that of
${}^u\Ind^{G}_{P}(\xi_{2\delta} )$.
(We can determine good orbits applying a combinatorial algorithm
described in [Trapa 2001].)
We denote by $\gca$ the set of the good open $G$-orbits on $X$.
We have:

{\bf Theorem E} \,\,\, (Corollary 2.10.3)
\[
{\rm Socle}({}^u\Ind_P^G(\xi_{2\delta}))=\bigoplus_{\hol\in\gca}\ana_\hol.
\]
Here, ${\rm Socle}({}^u\Ind_P^G(\xi_{2\delta}))$ is the maximal
semisimple submodule of ${}^u\Ind_P^G(\xi_{2\delta})$.
In this case, ${\rm Socle}({}^u\Ind_P^G(\xi_{2\delta}))$ is the image of
an intertwining operator induced from a homomorphism between
generalized Verma modules.

I thank Hiroyuki Ochiai for his informing me of typographical mistakes in earlier version.


\setcounter{section}{1}
\setcounter{subsection}{0}

\section*{\S\,\, 1.\,\,\,\, General connected complex reductive
groups and $\GL(n,\cpx)$}

\subsection{Notations}

Let $G$ be a connected complex reductive linear Lie group.
Let $\gggg$ be the Lie algebra of $G$, $U(\gggg)$ the universal enveloping algebra of $\gggg$, $\hhh$ a Cartan subalgebra of $\gggg$, and $\Delta$ the root system with respect to $(\gggg,\hhh)$.
Let $W$ be the Weyl group of the pair $(\gggg, \hhh)$ and 
denote by $w_0$ the longest element of $W$.
We fix a positive system $\Delta^+$ of $\Delta$ and denote by $\Pi$ the
corresponding basis of $\Delta$.
We fix a triangular decomposition $\gggg=\bnn\oplus\hhh\oplus\nnn$ such
that 
$\nnn = \sum_{\alpha\in\Delta^+}\gggg_\alpha$ and 
$\bnn = \sum_{-\alpha\in\Delta^+}\gggg_\alpha$.
Here, $\gggg_\alpha$ is the root space with respect to $\alpha$.
We consider a Borel subalgebra
$\bbb = \hhh+\nnn$.

We denote by $\gggg_0$ the normal real form of $\gggg$ which is compatible with the above decomposition and denote by $X\leadsto\bar{X}$ the complex conjugation with respect to $\gggg_0$.
Then there is an anti-automorphism, so-called Chevalley anti-automorphism,  $X\leadsto{}^tX$ of $\gggg$ which satisfies the following (1)-(3).

(1) \,\,\,\, ${}^t\gggg_0 =\gggg_0$.

(2) \,\,\,\, ${}^t\nnn =\bnn$,  ${}^t\bnn =\nnn$.

(3) \,\,\,\, ${}^tX = X$ \,\, $(X\in\hhh)$.

We extend $X\leadsto{}^tX$ to an anti-automorphism of $U(\gggg)$.

We define a homomorphism of a real Lie algebra $\gggg \rightarrow \gggg\times\gggg$ by $X\leadsto (X,\bar{X})$ for $X\in\gggg$.
Then the image of this homomorphism is a real form of $\gggg\times\gggg$.
Hence, we can regard $\gggg\times\gggg$ as the complexification $\ggc$ of $\gggg$.
$\kkc = \{(X,-^tX)\mid X\in\gggg\}$ is identified with the complexification of a compact form of $\gggg$.
$\kkc$ is also identified with $\gggg$ by an isomorphism $X\leadsto
(X,-{}^tX)$ as complex Lie algebras.
We may regard $G\times G$ as a complexification $G_\cpx$ of $G$.
Let $P$ be a standard parabolic subgroup of $G$; namely the Lie algebra $\ppp$ of $P$ contains $\bbb$.
We denote by $M$ (resp.\ $\mmm$) the Levi part of $P$ (resp.\ $\ppp$)
stable under the Chevalley anti-automorphism.
We denote by $\nnn_\ppp$ the nilradical of $\ppp$.
Let $S$ be the subset of $\Pi$ corresponding to $\ppp$; namely $S$ is
the basis of the root system with respect to $(\mmm, \hhh)$.
Put $S^\prime=\{-w_0\alpha|\alpha\in S\}$.
We denote by $\ppp^\prime$ the standard parabolic subalgebra
corresponding to $S^\prime$.
We denote by $w_\ppp$ the longest element of the Weyl
group with respect to $(\mmm, \hhh)$.

\subsection{Formulation of the problem and results}

Under the above identification $G_\cpx\cong G\times G$, the
complexification $P_\cpx$ of $P$ is identified with a subgroup $P\times P$ of $G\times G$.
The complex generalized flag variety $X=G_\cpx/P_\cpx$ is identified
with $G/P\times G/P$.
$X$ can be regarded as the set of parabolic subalgebra of $\gggg$ which is
$\Ad(G_\cpx)$-conjugate to $Lie(P_\cpx)=\ppp_\cpx\cong \ppp\times\ppp$.

We easily have:

\begin{prop}
$X$ has a unique $G$-orbit (say $\hol_0$).
The following five  conditions are equivalent to each other.

(1) \,\,\, $\hol_0$ contains a $\theta$-stable parabolic
subalgebra of $\gggg_\cpx$ .

(2) \,\,\, $\ppp$ and the parabolic subalgebra opposite to $\ppp$ are $G$-conjugate.

(3) \,\,\, $\ppp=\ppp^\prime$.

(4) \,\,\, $S$ is stable under the action of $-w_0$.

(5) \,\,\, $w_0 w_\ppp=w_\ppp w_0$. (Namely $w_0w_\ppp$ is an involution.)
\end{prop}

For $X\in \mmm$, we define 
$\delta_\ppp(X)=\frac{1}{2}\tr\left(\ad_\gggg(X)|_{\nnn_\ppp}\right)$.

Then, $\delta_\ppp$ is a one-dimensional representation of $\mmm$.
We see that $2n\delta_\ppp$ lifts to a holomorphic group homomorphism $\xi_{2n\delta_\ppp}
: M\rightarrow\cpx^\times$ for all $n\in\itg$.
Identifying $M_\cpx\cong M\times M$, for $n,m\in 2\itg$ we can regard
$\chi_{n,m}=\xi_{n\delta_\ppp}\boxtimes\xi_{m\delta_\ppp}$ as a real analytic
group homomorphism $M\rightarrow\cpx^\times$.
Taking account of  the natural projection $P\rightarrow M$, we regard
$\chi_{n,m}$ as a real analytic
one-dimensional representation of $P$.
So, we can consider an unnormalized degenerate principal series representation:
\[
 {}^u\Ind^G_P(\chi_{n,m})
=\{f\in C^\infty(G)|f(gp)=\chi_{n,m}(p)^{-1}f(g) \,\,\,\, (g\in G, p\in P)\}_{\mbox{$\kkc$-finite}}.
\]

If $\ppp=\bbb$, then
$\delta_\bbb=\rho=\frac{1}{2}\sum_{\alpha\in\Delta^+}\alpha$.
For $t\in\rel$, we define an unnormalized generalized Verma module:
\[
  {}^uM_\ppp(t)=U(\gggg)\otimes_{U(\ppp)}\cpx_{t\delta_\ppp}.
\]
Here, $\cpx_{t\delta_\ppp}$ means the one-dimensional module of $\ppp$
such that $\mmm$ acts on it by $t\delta_\ppp$.
If $\ppp=\bbb$, we simply write ${}^uM(t)={}^uM_\bbb(t)$.
We denote by $I_\ppp(t)$ the annihilator of $M_\ppp(t)$ in $U(\gggg)$.
From [Vogan 1984], $M_\ppp(t)$ is irreducible for $t\leqslant -1$.
Therefore $I_\ppp(t)$ is a primitive ideal for $t\leqslant -1$.

Under the identification $\ggc\cong\gggg\times\gggg$, we identify
$U(\ggc)\cong U(\gggg)\otimes U(\gggg)$.
We denote by $u \rightsquigarrow \check{u} \,\,\,  (u\in U(\gggg)\,)$
the anti-automorphism of $U(\gggg)$ generated by $X \leadsto
-X$ $(X\in\gggg)$.
If $V$ is a $U(\gggg)$-bimodule, we can regard $V$ as a
$U(\ggc)$-module as follows.
\[
  (u_1\otimes u_2)v = {}^t\check{u}_1v\check{u}_2 \,\,\,\,\,\,\,\,(u_1,u_2\in
  U(\gggg),v\in V\,).
\]
In particular, we may regard $U(\gggg)/I_\ppp(t)$ as a $U(\ggc)$-module.

We quote:
\begin{thm}\,\,$\mbox{([Conze-Berline, Duflo 1977] 2.12 and 6.3)}$

\mbox{}
If $n$ is an even integer such that $n\geqslant 2$, then ${}^u\Ind^G_P(\chi_{n,n})$ is isomorphic to
 $U(\gggg)/I_\ppp(-n)$ as $U(\ggc)$-modules.

\end{thm}

Together with [Borho-Kraft 1976] 3.6 and the primtivity of
$I_\ppp(-n)$ ($n\geqslant 0$), we have:

\begin{cor}
If $n\geqslant 2$, then ${}^u\Ind^G_P(\chi_{n,n})$ has a unique
 irreducible submodule (say $Y_n$).
The Gelfand-Kirillov dimension of $Y_n$ equals that of
 ${}^u\Ind^G_P(\chi_{n,n})$.
The Gelfand-Kirillov dimension of ${}^u\Ind^G_P(\chi_{n,n})/Y_n$ is
 strictly smaller than that of $Y_n$.
\end{cor}

[Borho-Jantzen 1977] 4.10 Corollar (also see
[Jantzen 1983] 15.27 Corollar) implies
$I_\ppp(-n)=I_{\ppp^\prime}(n-2)$.

Since ${}^u\Ind^G_P(\chi_{2,2})$ is the induced module associated with
the canonical bundle, we only consider the case of $n=2$.
For simplicity, we put $I=I_{\ppp}(-2)=I_{\ppp^\prime}(0)$.
Let $M$ and $N$ be $U(\gggg)$-modules.
$\Hom_\cpx(M,N)$ has natural structure of a $U(\gggg)$-bimodule.
We introduce a $U(\ggc)$-module structure on $\Hom_\cpx(M,N)$.
We denote the $\kkc$-finite part of $\Hom(M,N)$ by $L(M,N)$.
$V\rightsquigarrow L({}^uM(0),V)$ defines a category equivalence
between a category of highest weight modules and a category of
Harish-Chandra $(\ggc,\kkc)$-modules (Bernstein-Gelfand-Joseph-Enright
cf.\ [Bernstein-Gelfand 1980], [Joseph 1979]).
Via this category equivalence, the 2-sided ideals of $U(\gggg)/I$ correspond to the
$U(\gggg)$-submodules of ${}^uM(0)/I{}^uM(0)$.
In particular $Y_2$ corresponds to the irreducible highest weight module
 with the highest weight $\tau\rho-\rho$ ([Joseph 1979], also see
[Joseph 1983] p43). Here, $\tau\in W$ is the Duflo involution associated to the primitive
ideal $I$.
For integral weights $\mu,\nu\in\hhh^\ast$, we denote by $V(\nu,\mu)$ the Langlands (Zhelobenko) subquotient of
${}^u\Ind^G_B(\xi_{\nu+\rho}\boxtimes\xi_{\mu+\rho})$.
From [Joseph 1979] 4.5, we have:
\begin{thm}
The unique irreducible submodule $Y_2$ of ${}^u\Ind_P^G(\chi_{2,2})$ is
 isomophic to $V(-\tau\rho,-\rho)=V(w_0\tau w_0\rho,\rho)$.
\end{thm}

Next, under the equivalent conditions in Proposition 1.2.1, we compare $Y_2$ with
$\ana_{\hol_0}=H^{\dim\nnn_\ppp}(\hol_0,\lca)_{\mbox{$\kkc$-finite}}$.
First, we remark:

\begin{thm} \,\,\,\, ([Enright 1979], [Vogan-Zuckerman 1984])

\mbox{} Assume that $\ppp$ satisfies the equivalent conditions in Proposition 1.2.1.
Then we have
\[
\ana_{\hol_0} \cong 
{}^u\Ind^{G}_{P}(\chi_{2,0})
\cong
{}^u\Ind^{G}_{P}(\chi_{0,2}).
\]
\end{thm}

Since $2\delta_\ppp=\rho-w_0w_\ppp\rho$, we see $\ana_{\hol_0}$ is the
Langlands subquotient $V(-w_0w_\ppp\rho,-\rho)$ of
${}^u\Ind_B^G(\xi_{-w_0 w_\ppp\rho+\rho}\boxtimes\xi_{-\rho+\rho})$.

Since ${}^uM_\ppp(-2)$ is the irreducible highest weight module
of the highest weight $-2\delta_\ppp=w_0w_\ppp\rho-\rho$, we have:
\begin{thm}
\mbox{} Assume that $\ppp$ satisfies the equivalent conditions in Proposition 1.2.1.
Then $\ana_{\hol_0}$ is isomorphic to the unique submodule of
${}^u\Ind_P^G(\chi_{2,2})$ if and only if $w_0w_\ppp$ is a Duflo
involution in $W$.
\end{thm} 

In case that $w_0 w_\ppp$ is a Duflo involution, the embedding
 of $\ana_{\hol_0}$ into ${}^u\Ind^G_P(\chi_{2,2})$ can be regarded as intertwining operators between
degenerate principal series representation:
\begin{align*}
{}^u\Ind^{G}_{P}(\chi_{2,0}) & \hookrightarrow {}^u\Ind^G_P(\chi_{2,2}),
 \\
{}^u\Ind^{G}_{P}(\chi_{0,2}) & \hookrightarrow {}^u\Ind^G_P(\chi_{2,2}).
\end{align*}

In fact, the following result holds:

\begin{prop} \,\,\, ([Matumoto 1993])

Let $t$ be a non-negative even integer.
Then we have
\[
{}^uM_\ppp(-t-2)\hookrightarrow {}^uM_\ppp(t)
\]
if and only if $w_0 w_\ppp$ is a Duflo involution in $W$.
\end{prop}

Taking account of the isomorphism ${}^u\Ind^{G}_{P}(\chi_{m,n})\cong ({}^uM_\ppp(-m)\otimes {}^uM_\ppp(-n))^\ast_{\mbox{$\kkc$-finite}}$, the above intertwining operator
are induced from the embedding of the generalized Verma module.

If $W$ is a Weyl group of the type A, each involution in $W$ is a Duflo
involution ([Duflo 1977]).
Hence, we have:

\begin{cor}
If $G=\GL(n,\cpx)$ and $\ppp$ satisfies the equivalent conditions in
 Proposition 1.1, the socle of ${}^u\Ind^G_P(\chi_{2,2})$ is isomorphic
 to $\ana_{\hol_0}$.
\end{cor}


\setcounter{section}{2}
\section*{\S\,\, 2.\,\,\,\, $\U(m,n)$}
\setcounter{subsection}{0}

\subsection{{Root systems}}

Let $m$ and $n$ be non-negative integers.
We fix a maximal compact subgroup $K(m,n)$ of $G=\U(m,n)$, which is isomorphic
to $\U(m)\times \U(n)$.
We denote by $\theta$ the Cartan involution corresponding to $K(m,n)$.
We fix a $\theta$-stable maximally split Cartan subgroup $\sH(m,n)$ of
$\U(m,n)$.
We denote by $\shh(m,n)$ the complexified Lie algebra of $\sH(m,n)$.
We remark that all the Cartan subgroup of $\U(m,n)$ is connected.
We denote by $\gggg(m,n)$ (resp.\ $\kkk(m,n)$)the complexified Lie
algebra of $\U(m,n)$ (resp.\ $K(m,n)$).  
So, we have $\gggg(m,n)\cong\gl(m+n,\cpx)$ and
$\kkk(m,n)\cong\gl(m,\cpx)\times\gl(n,\cpx)$.

We denote by $K_\cpx(m,n)$ the complexification
of $K(m,n)$. So, $K_\cpx(m,n)$ is isomorphic to
$\GL(m,\cpx)\times\GL(n,\cpx)$.
We identify the complexification
of $\U(m,n)$ with $\GL(m+n,\cpx)$. 
First, we consider the root system $\Delta(m,n)$ with respect to $\Delta(\gggg(m,n),\shh(m,n))$.
Then we can choose an orthonormal basis $e_1,...,e_{m+n}$ of $\shh(m,n)^\ast$ such that 
\[ \Delta(m,n) = \{ e_i-e_j\mid 1\leqslant i,j\leqslant m+n,  i\neq
j\}.\]
We choose the above $e_1,...,e_{m+n}$ so that
$\theta(e_{i})=e_{m+n-i+1}$ for all $1\leqslant i\leqslant \min\{m,n\}$ and
$\theta(e_{i})=e_{i}$ for all $\min\{m,n\}< i< m+n-\min\{m,n\}$.

We fix a simple system $\Pi(m,n)$ for $\Delta(\gggg(m,n),\shh(m,n))$;
$\Pi=\{e_1-e_2,...,e_{m+n-1}-e_{m+n} \}$.
We denote by $\Delta^+(m,n)$ the corresponding positive system.

\subsection{Standard parabolic subgroups}

Fix  $\theta$, $\sH(m,n)$, etc. as in 2.1.

We also consider the particular orthonormal basis $e_1,...,e_{m+n}$ of
$\shh^\ast$ defined in 2.1.
We denote by $E_1,....,E_{m+n}$ the basis of $\shh(m,n)$ dual to $e_1,...,e_{m+n}$.

Let $\kappa=(k_1,...,k_s)$ be a finite sequence of positive integers
such that 
\[
 k=k_1+\cdots+k_s\leqslant \min\{m,n\}
.
\]
We put $k^\ast_i= k_1+\cdots+ k_i$ for $1\leqslant i\leqslant s$ and
$k^\ast_0=0$.
Hence $k^\ast_s=k$.
Put $m^\prime=m-k$ and
$n^\prime=n-k$.

We put
\[
 A_i=\sum_{j=1}^{2k_i}\left(E_{k^\ast_{i-1}+j}-E_{m+n-(k^\ast_{i-1}+j+1)}\right) \,\,\,\,\, (1\leqslant
 i\leqslant s)),¡¡
\]

Then we have $\theta(A_i)=-A_i$ for $1\leqslant
 i\leqslant s$.
We denote by $\aaa_\kappa(m,n)$  the Lie subalgebra of
$\shh(m,n)$ spanned by $\{A_i\mid 1\leqslant
 i\leqslant s\}$.

We define a subset $S(\kappa)$ of $\Pi(m,n)$ as follows.
\[ S(\kappa)=
\Pi-\{e_{k^\ast_i}-e_{k^\ast_i+1},e_{m+n-k^\ast_i}-e_{m+n-k^\ast_i+1} | 1\leqslant i\leqslant s\}.\]

We denote by $M_\kappa(m,n)$ (resp.\ $\mmm_\kappa(m,n)$) the standard Levi
subgroup (resp.\ subalgebra) of $\U(m,n)$ (resp.\ $\gggg(m,n)$) corresponding to
$S(\kappa)$.
Namely $M_\kappa(m,n)$ is the centralizer of $\aaa_\kappa(m,n)$ in $\U(m,n)$.

We denote by $P_\kappa(m,n)$ a parabolic subgroup of $\U(m,n)$ whose $\theta$-stable Levi part is $M_\kappa(m,n)$.
We choose $P_\kappa(m,n)$ so that the roots in $\Delta(m,n)$ whose root spaces are
contained in the complexified Lie algebra of the nilradical of
$P_\kappa(m,n)$ are all in $\Delta^+(m,n)$.
We denote by $N_\kappa(m,n)$ the nilradical of $P_\kappa(m,n)$.

Formally, we denote by $\U(0,0)$ the trivial group $\{1\}$ and we denote by $\GL(\kappa,\cpx)$ a product group $\GL(k_1,\cpx)\times\cdots\times \GL(k_s,\cpx)$.
Then, we have
\[ M_\kappa(m,n)\cong 
\GL(\kappa,\cpx)\times \U(m^\prime,n^\prime).\]
$\GL(\kappa,\cpx)$ and $\U(m^\prime,n^\prime)$ can be identified with a
subgroups of $M_\kappa(m,n)$ up to their automorphisms.
The Cartan involution $\theta$ induces  Cartan involutions on
$M_\kappa(m,n)$, $\GL(\kappa,\cpx)$, $\U(m^\prime, n^\prime)$ and we denote them by the same letter $\theta$.

We denote by $\ppp_\kappa(m,n)$, $\mmm_\kappa(m,n)$, and $\nnn_\kappa(m,n)$ the complexified Lie
algebras of $P_\kappa(m,n)$, $M_\kappa(m,n)$, and $N_\kappa(m,n)$,
respectively.

\subsection{Parabolic inductions}

For a  general reductive Lie group $G$ and its parabolic subgroup $P$, we
define an unnormalized induction as follows.
We fix a maximal compact subgroup $K$ and the corresponding Cartan
involution $\theta$ such that $P$ has a $\theta$-stable Levi part $L$.
We denote by $\gggg$, $\kkk$, and $\llll$ the complexified Lie algebras
of $G$, $K$,and $L$, respectively.
We denote by $U$ the nilradical of $P$.
Let $Z$ be a Harish-Chanda $(\llll, L\cap K)$-module and $(\pi, H)$ be a Hilbert
globalization of $Z$.
${}^u\Ind_{P}^G(Z)$  is the $K$-finite part of
\[\{f\in C^\infty(G)\otimes H\mid f(g\ell n)=\pi(\ell^{-1})f(g) \,\,\,\,\, (g\in G, \ell\in L, n\in U)\}.\]

Let $m$ and $n$ be non-negative integers and Let $k$ be a positive
integer such that $k\leqslant \min\{m,n\}$.
Let $\mu$ and $\nu$ be complex numbers such that $\mu+\nu\in\itg$.
We define a one-dimensional representation $(\eta^k_{\mu,\nu},
\cpx^k_{\mu,\nu})$ of $\GL(k,\cpx)$ as follows.
\[
 \eta^k_{\mu,\nu}(g)=\det(g)^\mu\overline{\det(g)}^{-\nu}
 \,\,\,\,\,\, (g\in \GL(k,\cpx)).
\]
When we identify $M_{(k)}(m,n)$ with $GL(k,\cpx)\times U(m-k,n-k)$,
there is ambiguity coming from automorphisms of $GL(k,\cpx)$ and 
$U(m-k,n-k)$.
We choose the identification so that the differential of the restriction of
$\cpx^k_{\mu,\nu}\boxtimes\cpx_{trivial}$ to $\sH(m,n)$ is
$\mu\sum_{i=1}^k e_i +\nu\sum_{j=1}^k e_{m+n-j+1}\in\shh(m,n)^\ast$.
For a Harish-Chandra $(\gggg(m-k,n-k),K(m-k, n-k))$-module $Z$, we define
\[
 {}^uI^k[\mu,\nu](Z)={}^u\Ind_{P_{(k)}(m,n)}^{U(m,n)}(\cpx^k_{\mu,\nu}\boxtimes
 Z).
\] 
Taking account of ``$\rho$-shift'', we define the normalized parabolic
induction by
${}^nI^k[\mu,\nu](Z)={}^uI^k[\mu+\frac{m^\prime+n^\prime+k}{2},\nu+\frac{m^\prime-n^\prime+k}{2}](Z)$.

\subsection{$\theta$-stable parabolic subalgebras}

In this article, we denote by $\nat$  the set of the non-negative integers.

We discuss $\theta$-stable parabolic subalgebras with respect to $\U(m,n)$
(cf.\ [Vogan 1996] Example 4.5).

Let $\ell$ be a positive integer. Put
\[\dbP_\ell(m,n)=\left\{((m_1,...,m_\ell),(n_1,...,n_\ell))\in \nat^{\ell}\times\nat^{\ell}\mid \sum_{i=1}^\ell m_i=m, \sum_{i=1}^\ell n_i=n, m_j+n_j>0 \;\; \mbox{for all $1\leq j\leq \ell$}\right\},\] 
We also put
\[\dbP(m,n)=\bigcup_{\ell> 0}\dbP_\ell(m,n),\]
\[\dbP(0,0)=\dbP_0(0,0)=\{((\emptyset),(\emptyset))\}.\]
If $(\bpm,\bpn)\in\dbP(m,n)$ satisfies $(\bpm,\bpn)\in\dbP_\ell(m,n)$,
we call $\ell$ the length of $(\bpm,\bpn)$. 
For $(\bpm,\bpn)\in\dbP(m,n)$, we define
\[
  I_{(\bpm,\bpn)}=\diag(I_{m_1},-I_{n_1},...,I_{m_\ell},-I_{n_\ell})
\]
Here, for a positive integer $I_k$ means the $k\times k$ identity
matrix.
If $k=0$, we regard $I_0$ as nothing. 
Then we can consider the following realization of $U(m,n)$.
\[
 \U(m,n)=\left\{g\in\GL(m+n,\cpx)| {}^t\bar{g}I_{(\bpm,\bpn)}g=I_{(\bpm,\bpn)}\right\}.
\]
Let $\theta$ be the Cartan involution given by the conjugation by $I_{(\bpm,\bpn)}$.
In this realization, we denote by $\qqq(\bpm,\bpn)$ the
block-upper-triangular parabolic subalgebra of
$\gggg(m,n)=\gl(m+n,\cpx)$ with blocks of sizes
$p_1+q_1,....,p_\ell+q_\ell$ along the diagonal.
Then, $\qqq(\bpm,\bpn)$ is a $\theta$-stable parabolic subalgebra.
The corresponding Levi subgroup $\U(\bpm,\bpn)$ consists of diagonal
blocks.
\[
\U(\bpm,\bpn)\cong \U(m_1,n_1)\times\cdots\times\U(m_\ell,n_\ell). 
\]
We denote by $\gggg(\bpm,\bpn)$ the complexified Lie algebra of $\U(\bpm,\bpn)$.
We denote by $\vvv(\bpm,\bpn)$ the nilradical of $\qqq(\bpm,\bpn)$.

Via the above construction of $\qqq(\bpm,\bpn)$, $K(m,n)$-conjugate class of $\theta$-stable parabolic subalgebras with
respect to $\U(m,n)$  is classified by $\dbP(m,n)$.

\subsection{{Cohomological inductions}}

First, we discuss the case of $\theta$-stable maximal parabolic
subalgbras.

We fix non-negative integers $m$ and $n$ and
$((m^\prime,p),(n^\prime,q))\in\dbP_2(m,n)$. Here, $m^\prime=m-p$ and  $n^\prime=n-q$.
For an integer $h$, we define a unitary character
$\eta_{h}^{p,q} : \U(p, q) \rightarrow \cpx^\times$ by
\[
 \eta_{h}^{p,q}(g)=\det(g)^{h}
 \,\,\,\,\,\,\, (g\in\U(p,q)).
\] 
We denote by $\cpx^{p,q}_h$ the representation space of $\eta_{h}^{p,q}$.
For a Harish-Chandra $(\gggg(m-p,n-q),K(m-p,n-q))$-module $Z$, we define
\[
 R^{p,q}[h](Z)=\left({}^u\rca_{\qqq((m^\prime,p),(n^\prime,q)),K(m^\prime,n^\prime)\times
 K(p,q)}^{\gggg(m,n),K(m,n)}\right)^{pm^\prime+qn^\prime}(Z\boxtimes\cpx^{p,q}_{h+m+n})\otimes\cpx^{m,n}_{p+q}.
\]
Here,
$\left({}^u\rca_{\qqq((m^\prime,p),(n^\prime,q)),K(m^\prime,n^\prime)\times
K(p,q)}^{\gggg(m,n),K(m,n)}\right)^{pm^\prime+qn^\prime}$
is an unnormalized cohomological induction defined in [Knapp-Vogan 1995]
p676.
If
$\left(\rca_{\qqq((m^\prime,p),(n^\prime,q)),K(m^\prime,n^\prime)\times K(p,q)}^{\gggg(m,n),K(m,n)}\right)^{pm^\prime+qn^\prime}$
is a cohomological induction functor defined in [Knapp-Vogan 1995] p328 or
[Vogan 1982], then we
have
\[
R^{p,q}[h](Z)=
 \left(\rca_{\qqq((m^\prime,p),(n^\prime,q)),K(m^\prime,n^\prime)\times K(p,q)}^{\gggg(m,n),K(m,n)}\right)^{pm^\prime+qn^\prime}(Z\boxtimes\cpx^{p,q}_{h}).
\]

Let $(\bpm,\bpn)=((m_1,...,m_\ell),(n_1,...,n_\ell))\in\dbP_\ell(m,n)$
and $\bph=(h_1,....,h_\ell)\in\itg^\ell$.
We define
\[
 \aca_{(\bpm,\bpn)}[\bph]=R^{m_\ell,n_\ell}[h_\ell](R^{m_{\ell-1},n_{\ell-1}}[h_{\ell-1}](\cdots(R^{m_{2},n_{2}}[h_{2}](\cpx^{m_1,n_1}_{h_1}))\cdots)).
\]
$\aca_{(\bpm,\bpn)}[\bph]$ is nothing but a derived functor module for
$\qqq(\bpm,\bpn)$.

We quote:

\begin{thm} \,\, \mbox{([Vogan 1984], [Vogan 1988], [Trapa 2001] Lemma 3.5)}

(1) \,\,  $\aca_{(\bpm,\bpn)}[\bph]$ is a derived functor module in the good
 range for $\qqq(\bpm,\bpn)$ if and only if $h_i\geqslant h_{i+1}$ for
 all $1\leqslant i<\ell$.
In this case, $\aca_{(\bpm,\bpn)}[\bph]$ is non-zero and irreducible.

(2) \,\, $\aca_{(\bpm,\bpn)}[\bph]$ is a derived functor module in the
 weakly fair
 range for $\qqq(\bpm,\bpn)$ if and only if 
\[
 h_i-h_{i+1}\geqslant -\frac{m_i+n_i+m_{i+1}+n_{i+1}}{2}
\,\,\,\,\, (1\leqslant i<\ell).\]

(3) \,\, $\aca_{(\bpm,\bpn)}[\bph]$ is a derived functor module in the
 mediocre
 range for $\qqq(\bpm,\bpn)$ if and only if 
\[
 h_i-h_{i+1}\geqslant -\max\{m_i+n_i, m_j+n_j\}-\sum_{i<k<j}(m_k+n_k)
\,\,\,\,\, (1\leqslant i<j\leqslant\ell).\]
In this case $\aca_{(\bpm,\bpn)}[\bph]$ is zero or irreducible.

The mediocre condition is weaker than the weakly fair condition.
The weakly fair condition is weaker than the good condition.
The mediocre condition implies the unitarity of $\aca_{(\bpm,\bpn)}[\bph]$.

\end{thm}

\subsection{Associated varieties}

It is well-known that the $\GL(m+n,\cpx)$-nilpotent orbits in
$\gggg(m,n)=\gl(m+n,\cpx)$ are parametrized by the  Young diagrams with
$m+n$ boxes. (Cf.\ [Collingwood-McGovern 1993])  

Fix a sequence of positive integers $\bpc=(c_1,...,c_\ell)$ such that
$c_1+\cdots+c_\ell=m+n$.
Let $\tau\in\gS_\ell$ be a permutation such that $c_{\tau(1)}\geqslant
c_{\tau(2)}\geqslant\cdots\geqslant c_{\tau(\ell)}$.
Denote by $Y(\bpc)$ the Young diagram consisting of $\ell$ columns such
that the length of the $i$-th column of
$Y(\bpc)$ is $c_{\tau(i)}$.
We denote by $\hol(\bfc)$ the nilpotent $\GL(m+n,\cpx)$-orbit in
$\gggg(m,n)$ corresponding to $Y(\bpc)$.

For $(\bpm,\bpn)=((m_1,...,m_\ell),(n_1,...n_\ell))\in\dbP_\ell(m,n)$,
we put $\bpm+\bpn=(m_1+n_1,...,m_\ell+n_\ell)$.
For $(\bpm,\bpn)=((m_1,...,m_\ell),(n_1,...n_\ell))\in\dbP_\ell(m,n)$, the
Richardson orbit with respect to $\qqq(\bpm,\bpn)$ coincides with $\hol(\bpm+\bpn)$. (Kraft, Ozeki, Wakimoto,
see [Collingwood-McGovern 1993] Theorem 7.2.3 ) 
Namely, $\hol(\bpm+\bpn)$ is a unique nilpotent $\GL(m+n,\cpx)$-orbit in
$\gggg(m,n)$ such that $\hol(\bpm,\bpn)\cap\vvv(\bpm,\bpn)$ is open in
$\vvv(\bpm,\bpn)$.
(Here, $\vvv(\bpm,\bpn)$ is the niradical of $\qqq(\bpm,\bpn$.)
It is known that $\dim\hol(\bpm+\bpn)=2\dim\vvv(\bpm,\bpn)$.

Let $\gggg(m,n)=\kkk(m,n)\oplus\sss(m,n)$ be the
complexified Cartan decomposition with respect to $\theta$.

The nilpotent $K_\cpx(m,n)$-orbits in $\sss(m,n)$ are parametrized by (the equivalence
classes of) the signed Young
diagrams of signature $(m,n)$. (See [Collingwood-McGovern 1993])
Here, a signed Young diagram is a Young diagram in which every box is
labeled with $+$ or $-$ sign in such a way that signs alternate across
rows.
We call two signed Young diagrams equivalent, if one can be obtained
from the other by interchanging rows of equal length. 
The signature of a signed Young diagram is the ordered pair $(m,n)$,
where $m$ (resp.\ $n$) is the number of boxes labeled $+$ (resp.\ $-$).
For a signed Young diagram $T$, we denote by $|T|$ the Young diagram
obtained by earasing all $+$ and $-$ on $Y$.
We call $|T|$ the shape of $T$.
For a Young diagram consisting of $m+n$ boxes, we denote by
$\dbS_{m,n}(Y)$ the set of the equivalence classes of the signed Young diagrams of signature
$(m,n)$ whose shape is $Y$.

The above parametrization of the nilpotent $K_\cpx(m,n)$-orbits in $\sss(\bpm,\bpn)$ is compatible with the parametrization of
$\GL(m+n,\cpx)$-nilpotent orbits in
$\gggg(m,n)=\gl(m+n,\cpx)$ via Young diagrams.
Namely, for any  nilpotent $K_\cpx(m,n)$-orbit $\hol$ in $\sss(m,n)$
corresponding to a signed Young diagram $T$,
the nilpotent $\GL(n+m,\cpx)$-orbit $\Ad(\GL(m+n,\cpx))\hol$
corresponds to $|T|$.
For a nilpotent $\GL(m+n,\cpx)$-orbit $\hol$ in $\gggg(m,n)$
corresponding to a Young diagram $Y$, $\hol\cap\sss(m,n)$ is the union
of the nilpotent $K_\cpx(m,n)$-orbits corresponding to the elements of 
$\dbS_{m,n}(Y)$. 

For a Harish-Chandra $(\gggg(m,n), K(m,n))$-module $M$, we denote by
$\Dim(M)$ (resp.\ $\Ass(M)$) the Gelfand-Krillov dimension (resp.\ the
associated variety) of $M$. (cf.\ [Vogan 1978], [Vogan 1991])
For an irreducible Harish-Chandra $(\gggg(m,n), K(m,n))$-module $M$ with
an integral infinitesimal character, $\Ass(M)$ is the closure of a single
nilpotent $K_\cpx(m,n)$-orbit in $\sss(m,n)$. ([Barbasch-Vogan 1982])

\begin{lem} \,\,\,([Trapa 2001] Lemma 5.6)
\,\,\, \mbox{}

Let $Z$ be a Harish-Chandra $(\gggg(m,n), K(m,n))$-module such that
 $\Ass(Z)$ is the closure of a single nilpotent $K_\cpx(m,n)$-orbit  
 in $\sss(m,n)$ corresponding to a signed Young diagram $T$.
Let $h$ be a sufficiently large integer so that the ccohomological
 induction $R^{p,q}[h](Z)$ is in the good range.
Then $\Ass(R^{p,q}[h](Z))$ is the closure of a nilpotent
 $K_\cpx(m+p,n+q)$-orbit (say $\hol$) in $\sss(m+p,n+q)$.
Then, the signed Young diagram corresponding to $\hol$ is obtained by
 adding $p$ pluses and $q$ minuses, from top to bottom, to the row-ends
 of $T$ so that

 (1) \,\,\, at most one sign is added to each row-end,

 (2) \,\,\, the signs of the resulting diagram must alternate across
 rows, and

 (3) \,\,\, add signs  to the highest row possible.

(The resulting diagram may not have rows of decreasing length.
 If necessary, we can choose another signed Young diagram equivalent to $T$  in order
 to obtain a signed Young diagram as the result of the above procedure.)

\end{lem}

{\it Remark}  \,\,\,\, As explained in [Trapa 2001], we can calculate
the associated varieties of derived functor modules of $\U(m,n)$ using Lemma 2.6.1.

For $(\bpm,\bpn)\in\dbP(m,n)$, $\Dim(\aca_{\bpm,\bpn}[\bph])\leqslant\frac{1}{2}\dim\hol(\bpm,\bpn)$.
We call $(\bpm,\bpn)\in\dbP(m,n)$ normal, if
$\Dim(\aca_{\bpm,\bpn}[\bph])=\frac{1}{2}\dim\hol(\bpm,\bpn)$ holds in
the good range. It is kown that for a normal $(\bpm,\bpn)\in\dbP_\ell(m,n)$ and any mediocre $\bph\in\itg^\ell$, we have
$\Dim(\aca_{\bpm,\bpn}[\bph])=\frac{1}{2}\dim\hol(\bpm,\bpn)$ and
$\Ass(\aca_{\bpm,\bpn}[\bph])=\Ass(\aca_{\bpm,\bpn}[\bfh^\prime])$,
where $\bfh^\prime$ is any parameter in the good range. (cf.\ [Trapa 2001])
Fix a sequence of positive integers $\bpc=(c_1,...,c_\ell)$ such that
$c_1+\cdots+c_\ell=m+n$.
Put
\[  
\dbO(\bpc)=\{(\bpm,\bpn)\in\dbP_\ell(m,n)\mid
\mbox{$(\bpm,\bpn)$ is normal and $m_i+n_i=c_i$ for all $ 1\leqslant i\leqslant\ell$} \}
\]

Examining the procedure described in Lemma 2.6.1, we can easily have:

\begin{cor}
Let $\bpc$ be a sequence of positive integers as above.
The correspondence
 $(\bpm,\bpn)\rightsquigarrow\Ass(\aca_{(\bpm,\bpn)}[\bph])$ induces a
 bijection
\[
 \Phi_{\bpc} : \dbO(\bpc)\stackrel{\sim}{\rightarrow} \dbS_{m,n}(Y(\bpc)).
\]
\end{cor}

\subsection{The annihilator of a degenerate principal series representation}

Let $m$ and $n$ be non-negative integers.
Let $\kappa=(k_1,...,k_s)$ be a finite sequence of positive integers
such that 
$k=k_1+\cdots+k_s\leqslant \min\{m,n\}$.
Let $\bpu=(u_1,...,u_s)$ and $\bpv=(v_1,...,v_s)$ be sequences of complex
numbers such that $u_i+v_i\in\itg$ for all $1\leqslant i\leqslant s$.
Let $h\in\itg$.
We consider a degenerate principal series representation of $\U(m,n)$ as follows. 
\[
 {}^nI^\kappa_{m,n}[\bpu;h;\bpv]={}^nI^{k_1}[u_1,v_1]({}^nI^{k_2}[u_2,v_2](\cdots({}^nI^{k_s}[u_s,v_s](\cpx^{m-k,n-k}_h)\cdots)).
\]
If $m=n=k$, we regard $\cpx^{0,0}_h$ as the trivial representation of
the trivial group.
In this case ${}^nI^\kappa_{m,n}[\bpu;h;\bpv]$ does not depend on $h$.
We denote by  $\delta^\kappa_{m,n}$ an element of $\cpx^s$ whose 
$i$-th entry is $\frac{m+n-k^\ast_{i-1}}{2}$, where $k^\ast_i=k_1+\cdots+k_i$. 
We denote by $J(\bfu,h,\bfv)$ the annihilator of
${}^nI^\kappa_{m,n}[\bpu;h;\bpv]$ in the universal enveloping algebra $U(\gggg(m,n))$.

We define a normalized generalized Verma module as follows.
Let $\ppp$ be a parabolic subalgebra of $\gggg=\gggg(m+n)$.
Let $\ppp=\llll+\nnn$ be a Levi decomposition of $\ppp$.
We define a one-dimensional representation $-\rho_\ppp$ of $\llll$ by
$-\rho_\ppp(X)=\frac{1}{2}(\ad(X)|_{\nnn}) \,\,\,\,\, (X\in\llll)$.
For a one-dimensional representation $\xi$ of $\llll$ we extend
$\xi\otimes-\rho_\ppp$ to a one-dimensional representation of $\ppp$
as usual and define ${}^nM_\ppp(\xi)=U(\gggg)\otimes_{U(\ppp)}(\xi\otimes-\rho_\ppp)$.

We define a weight on $\shh(m,n)$ by 
\[
 \xi(\bfu,h,\bfv)=\sum_{i=1}^su_i\left(\sum_{j=1}^{k_i}e_{k^\ast_{i-1}+j}\right)+h\sum_{j=1}^{m+n-2k}e_{k+j}+\sum_{i=1}^sv_{s-i}\left(\sum_{j=1}^{k_i}e_{m+n-k+k_{i-1}^\ast+j}\right).
\]
Obviously, $\xi(\bfu,h,\bfv)$ can be extend to a one-dimensional
representation of $\llll$.
We denote by $\bpp_\kappa(m,n)$ the opposite parabolic subalgebra to $\ppp_\kappa(m,n)$.
${}^nI^\kappa_{m,n}[\bpu;h;\bpv]$ and ${}^nM_{\ppp_\kappa(m,n)}(-\xi(\bfu,h,\bfv))$ have a
perfect pairing. ${}^nM_{\bpp_\kappa(m,n)}(\xi(\bfu,h,\bfv))$ and
${}^nM_{\ppp_\kappa(m,n)}(-\xi(\bfu,h,\bfv))$ also have a perfect
pairing.

Hence we have:
\begin{lem}
$J(\bfu,h,\bfv)$ coincides with the annihilator of
 ${}^nM_{\bpp_\kappa(m,n)}(\xi(\bfu,h,\bfv))$.
\end{lem}

Let $\bfc=(c_1,...,c_\ell)$ be a sequence of non-negative integers such that
$c_1+\cdots+c_\ell=m+n$.
Put $c^\ast_i=c_1+\cdots+c_i$ and $c^\ast_{0}=0$.
We define a subset $S[\bfc]$ of $\Pi(m,n)$ as follows.
\[ S[\bfc]=
\Pi-\{e_{c^\ast_i}-e_{c^\ast_i+1} | 1\leqslant i\leqslant \ell\}.\]
We denote by $\ppp(\bfc)$ the standard parabolic subalgebra of
$\gggg(m,n)$ associated to $ S[\bfc]$.
Under appropriate realization of $\gggg(m,n)$ as $\gl(m+n,\cpx)$, $\ppp(\bfc)$
is the block-upper-triangular parabolic subalgebra of $\gl(m+n,\cpx)$
with blocks of sizes $c_1,...,c_\ell$ along the diagonal.
If $c_i=0$ for some $i$, ``a block of size $0$'' means nothing.
Namely, we simply neglect it.
We denote by $\bpp(\bfc)$ the parabolic subalgebra opposite to $\ppp(\bfc)$.
For a sequence $\bfh=(h_1,...,h_\ell)$ of complex numbers, we define a
weight on $\shh(m,n)$ as follows.
\[
 \xi(\bfh)=\sum_{i=1}^\ell h_i\left(\sum_{j=1}^{c_i}e_{c^\ast_{i-1}+j}\right).
\]

Hereafter, we consider a degenerate principal series representation
${}^nI^\kappa_{m,n}[\bpu;h;\bpv]$ with an integral infinitesimal
character.
Namely, we assume that
$\bfu-\delta^\kappa_{m,n},\bfv+\delta^\kappa_{m,n}\in\itg^s$.
We define $\bfh=(h_1,...,h_{2s+1})$ as follows.
\begin{align*}
 h_1& =h, \\
 h_{2i} &=u_{s-i+1} \,\,\,\, (1\leqslant i\leqslant s),\\
 h_{2i+1} &=v_i \,\,\,\, (1\leqslant i\leqslant s).
\end{align*}
We also define a sequence of positive integer $\bfc$ by
\begin{align*}
c_1 &=m+n-2k,\\
c_{2i} &= c_{2i+1}=k_{s-i+1} \,\,\,\,\, (1\leqslant i\leqslant s).
\end{align*}

We define
\[
 \Xi(\bfu,h,\bfv)=\{\tau\in{\mathfrak S}_{2s+1}\mid
 h_{\tau(1)}\geqslant\cdots\geqslant h_{\tau(2s+1)}\}.
\]
Obviously $\Xi(\bfu,h,\bfv)$ is non-empty.
For $\tau\in{\mathfrak S}_{2s+1}$, put
$\bfc^\tau=(c_{\tau(1)},...,c_{\tau(2s+1)})$ and $\bfh^\tau=(h_{\tau(1)},...,h_{\tau(2s+1)})$
From [Borho-Jantzen 1977] 4.10 Corollar and Lemma 2.7.1, we have:

\begin{prop}
If $\tau\in\Xi(\bfu,h,\bfv)$,
then $J(\bfu,h,\bfv)$ coincides with the annihilator of
 ${}^nM_{\bpp(\bfc^\tau)}(\xi(\bfh^\tau))$.
\end{prop}

From [Vogan 1984], we see that ${}^nM_{\bpp(\bfc^\tau)}(\xi(\bfh^\tau))$ is irreducible
for $\tau\in\Xi(\bfu,h,\bfv)$.
Hence we see that $J(\bfu,h,\bfv)$ is a primitive ideal.

\subsection{The associated variety of a degenerate principal series
  representation}

For a Harish-Chandra $(\gggg(m,n), K(m,n))$-module $Z$, we denote by
$[Z]$ the distribution character of $Z$.
Let $V$ and $M$ be Harish-Chandra $(\gggg(m,n), K(m,n))$-modules.
We write
$V\sim M$, if a virtual character $[V]-[M]$ is a linear
combination of  the distribution characters of irreducible Harish-Chandra
$(\gggg(m,n), K(m,n))$-modules whose Gelfand-Kirillov dimensions are strictly smaller
than $\Dim(V)$.

As in 2.7, we consider a degenerate principal series representation
${}^nI^\kappa_{m,n}[\bpu;h;\bpv]$ with an integral infinitesimal
character. We also define $\bfc$ and $\bfh$ from $(\bpu,h,\bpv)$ as in
2.7.
We have:

\begin{prop}
Assume that $\Xi(\bpu,h,\bpv)$ contains the identity element $e$.
Then, we have
\[
 {}^nI^\kappa_{m,n}[\bpu;h;\bpv]\sim\bigoplus_{(\bfm, \bfn)\in\dbO(\bfc)}\ana_{(\bfm,\bfn)}[\check{\bfh}]
\]
Here, $\check{\bfh}=(\check{h}_1,...,\check{h}_{2s+1})$ is defined by
 $\check{h}_i=h_i-\frac{m+n-c_i}{2}+c^\ast_{i-1}$ \,\,\,\,
 $(1\leqslant i\leqslant 2s+1)$.
\end{prop}

\proof
If $h_{2i}=h_{2i+1}$ for $1\leqslant i\leqslant s$, the irreducible
decomposition of ${}^nI^\kappa_{m,n}[\bpu;h;\bpv]$  is given in [Matumoto
1996] Theorem 3.3.1.
Each irreducible component is a derived functor module.
So, the statement of the proposition in this particular case is
obtained by examining the method of computation, via Lemma 2.6.1, of the associated
varieties of derived functor modules appearing the decomposition
formula.
The general case is obtained by applying the standard argument of the
translation principle into the weakly fair range. ([Vogan 1984], [Vogan 1988], [Vogan 1990])
\,\,\,\, $\blacksquare$

{\it Remark} \,\,\,\, We may prove Proposition 2.8.1 directly using
[Matumoto 2002] Theorem 2.2.3.

Remark that $\hol(\bfc)$ is the Richardson orbit with respect to
$\ppp_\kappa(m,n)$.
We denote by $\overline{\hol(\bfc)}$ its closure.

It is easy to see that $\Ass({}^nI^\kappa_{m,n}[\bpu;h;\bpv])$ does not depend on the
choice of the parameter $(\bpu;h;\bpv)$ such that
${}^nI^\kappa_{m,n}[\bpu;h;\bpv])$  has an integral infinitesimal
character.
Hence, we have:

\begin{cor}

For any ${}^nI^\kappa_{m,n}[\bpu;h;\bpv]$ with an integral infinitesimal
character, we have
\[
 \Ass({}^nI^\kappa_{m,n}[\bpu;h;\bpv])= \overline{\hol(\bfc)}\cap\sss(m,n)
\]

\end{cor}

From [Schmid-Vilonen 2000], we can identify the wave front sets of
derived functor modules.  The wave front set of the degenerate
principal series represntation is obtained by [Barbasch-Vogan 1980].
Hence, Corollary 2.8.2 is equivalent to the following result.

\begin{cor}
The real Lie algebra $Lie(N_\kappa(m,n))$ intersects all the real form
 of the Richardson $\GL(m+n,\cpx)$-orbit in $\gl(m+n,\cpx)$
 with respect to $\ppp_\kappa(m,n)$.
\end{cor}

{\it Remark } \,\,\,
The corresponding statement to Corollary 2.8.3 is
incorrect for a general semisimple Lie algebra.
For instance, the Jacobi parabolic subalgebra  of ${\sss}{\ppp}(2,\rel)$
gives  a counterexample; there are three real forms but the nilradical
intersects only one of them.

\subsection{Irreducible constituents of the maximal Gelfand-Kirillov
 dimension}

As in 2.7, we consider a degenerate principal series representation
${}^nI^\kappa_{m,n}[\bpu;h;\bpv]$ with an integral infinitesimal
character. We also define $\bfc$ and $\bfh$ from $(\bpu,h,\bpv)$ as in
2.7.

We have:

\begin{thm}
For $\tau\in\Xi(\bpu,h,\bpv)$, we have
\[
 {}^nI^\kappa_{m,n}[\bpu;h;\bpv]\sim\bigoplus_{(\bfm, \bfn)\in\dbO(\bfc^\tau)}\ana_{(\bfm,\bfn)}[\check{\bfh}^\tau]
\]
Here, $\check{\bfh}^\tau=(\check{h}_1^\tau,...,\check{h}_{2s+1}^\tau)$ is defined by
 $\check{h}_i^\tau=h_{\tau(i)}-\frac{m+n-c_{\tau(i)}}{2}+{c^\tau}^\ast_{i-1}$ \,\,\,\,
 $(1\leqslant i\leqslant 2s+1)$,
where ${c^\tau}^\ast_{i-1}=c_{\tau(1)}+\cdots+c_{\tau(i-1)}$ for
 $2\leqslant i\leqslant 2s+1$ and ${c^\tau}^\ast_{0}=0$.
\end{thm}

\proof

Applying the standard argument of the
translation principle into the weakly fair range ([Vogan 1984], [Vogan
1988], [Vogan 1990]), we may assume that $h_{\tau(1)}\gg\cdots\gg
h_{\tau(2s+1)}$.

From [Barbasch-Vogan 1983] (also see [Trapa 2001] section 6), we see

(1) \,\,\, An irreducible Harish-Chandra $(\gggg(m,n),K(m,n))$-module
    with an integral infinitesimal character  is
    determined by its annihilator and its associated variety.

(2) \,\,\, The associated variety of an irreducible Harish-Chandra $(\gggg(m,n),K(m,n))$-module
    with an integral infinitesimal character is irreducible.

Let $(\bfm, \bfn)\in\dbO(\bfc^\tau)$.
From [Vogan 1986] Proposition 16.8, we see the annihilator of
$\ana_{(\bfm,\bfn)}[\check{\bfh}^\tau]$ contains the annihilator of
 ${}^nM_{\bpp(\bfc^\tau)}(\xi(\bfh^\tau))$.
Since
$\Dim(\ana_{(\bfm,\bfn)}[\check{\bfh}^\tau])=\frac{1}{2}\dim\hol(\bfc)$,
Proposition 2.7.2 implies that the annihilator of
$\ana_{(\bfm,\bfn)}[\check{\bfh}^\tau]$ coincides with $J(\bfu,h,\bfv)$.
So, taking account of the above (1) and (2), we have
\[
 {}^nI^\kappa_{m,n}[\bpu;h;\bpv]\sim\bigoplus_{(\bfm,
 \bfn)\in\dbO(\bfc^\tau)}\ana_{(\bfm,\bfn)}[\check{\bfh}^\tau]^{\oplus n(\bfm,\bfn)}. 
\]
Here, $n(\bfm,\bfn)$ are positive integers.

Hence, we have only to show that:

($\Diamond$) \,\,\,\, The number of irreducible constituents
of the maximal Gelfand-Kirillov dimension in
${}^nI^\kappa_{m,n}[\bpu;h;\bpv]$ equals the number of signed Young
diagrams of the shape $Y(\bfc)$. 

Since the datum $(\bpu,h,\bpv)$ is recovered from $\bfh$, we write
$I[\bfh]={}^nI^\kappa_{m,n}[\bpu;h;\bpv]$, $J(\bfh)=J(\bpu,h,\bpv)$, and
 $\Xi[\bfh]=\Xi(\bpu,h,\bpv)$.
  
We regard the symmetric group ${\mathfrak S}_{2s+1}$ as a Coxeter group
such that a set of transpositions $\{(i,i+1)\mid 1\leq i<2s+1\}$ is the
simple reflections.
We denote by $\ell(\tau)$ the length of $\tau\in{\mathfrak S_{2s+1}}$.

We prove ($\Diamond$) by the induction on $\ell(\tau)$.
The case of $\ell(\tau)=0$ follows from Proposition 2.8.1.

So, we assume $\ell(\tau)>0$.

Then there exists some $1\leqslant i\leqslant 2s+1$ such that
the length of $\tau^\prime=\tau\circ(i,i+1)$ is strictly smaller than $\tau$.
Since we assumed that $h_{\tau(i)}\gg h_{\tau(i+1)}\gg h_{\tau(i+2)}$,
there is a positive integer $r_0$ such  that $h_{\tau(i+1)}\gg
h_{\tau(i)}-r_0\gg h_{\tau(i+2)}$.
For $r\in\nat$, we define $\bfh(r)=(\bfh(r)_1,...,\bfh(r)_{2s+1})$ by
$\bfh(r)_j=h_j$ for $j\neq\tau(i)$ and
$\bfh(r)_{\tau(i)}=h_{\tau(i)}-r$.
Then, we see  $\Xi(\bfh)=\{\tau\}$ and $\Xi(\bfh(r_0))=\{\tau^\prime\}$.
For $r\in\nat$, we define weights on $\shh(m,n)$ as follows.
\begin{align*}
\psi &= \sum_{\substack{1\leqslant j\leqslant 2s+1 \\ j\neq
 i,i+1}}\sum_{d=1}^{c_{\tau(j)}}\left(h_{\tau(j)}+\frac{c_{\tau(j)}-1}{2}-d\right)e_{{c^\tau}^\ast_{j-1}+d}\\
\psi(r)&=\psi+\sum_{d=1}^{c_{\tau(i)}}\left(h_{\tau(i)}-r+\frac{c_{\tau(i)}-1}{2}-d\right)e_{{c^\tau}^\ast_{i-1}+d}+\sum_{d=1}^{c_{\tau(i+1)}}\left(h_{\tau(i+1)}+\frac{c_{\tau(i+1)}-1}{2}-d\right)e_{{c^\tau}^\ast_{i}+d}\\
\psi^\prime(r)&=\psi+\sum_{d=1}^{c_{\tau(i+1)}}\left(h_{\tau(i+1)}+\frac{c_{\tau(i+1)}-1}{2}-d\right)e_{{c^\tau}^\ast_{i-1}+d}+\sum_{d=1}^{c_{\tau(i)}}\left(h_{\tau(i)}-r+\frac{c_{\tau(i)}-1}{2}-d\right)e_{{c^\tau}^\ast_{i-1}+c_{\tau(i+1)}+d}.
\end{align*}
Then, $\psi(r)$ and $\psi^\prime(r)$ are contained in the same Weyl
group orbit and represent the same infinitesimal character.
In fact, they are the infinitesimal character of $I[\bfh(r)]$.

For $\lambda,\mu\in\shh(m,n)^\ast$ such that $\lambda-\mu$ is integral, we
denote by $T^\mu_\lambda$ the translation functor.
Namely, if $Z$ is a Harish-Chandra module with a generalized infinitesimal character
$\lambda$, $T^\mu_\lambda(Z)=P_{\mu}(V_{\mu-\lambda}\otimes_\cpx Z)$.
Here, $V_{\mu-\lambda}$ is the irreducible finite-dimensional
representation with an extreme weight $\mu-\lambda$ and $P_\mu$ is the
projection functor to the generalized infinitesimal character $\mu$.

For non-negative integers $r$ and $d$ such that $r>d$, we put
\begin{align*}
\tca_d^r &=T^{\psi(r)}_{\psi(r-1)}\circ
 T^{\psi(r-1)}_{\psi(r-2)}\circ\cdots\circ T^{\psi(d+1)}_{\psi(d)}, \\
{\tca^\prime}_r^d &=T^{\psi^\prime(d)}_{\psi^\prime(d+1)}\circ
 T^{\psi^\prime(d+1)}_{\psi^\prime(d+2)}\circ\cdots\circ
 T^{\psi^\prime(r-1)}_{\psi^\prime(r)}.
\end{align*}

Using the idea of the proof of Trapa's result([Trapa 2001] Lemma 3.13) together with [Vogan
1982] Lemma 7.2.9 (b),  [Vogan 1988] Lemma 4.8, and Proposition 2.7.2
above, we also have:

\begin{lem} \,\,\,\, \mbox{}

(1) \,\,\,\, If $r\in\nat$ satisfies $r\leqslant h_{\tau(i)}-h_{\tau(i+1)}+\left|\frac{c_{\tau(i)}-c_{\tau(i+1)}}{2}\right|$, we have
$[\tca^r_0(I[\bfh])]= [I[\bfh(r)]]$.

(2) \,\,\,\, We fix $r_0\in\nat$ as above.
If $r\in\nat$ satisfies $r\geqslant h_{\tau(i)}-h_{\tau(i+1)}-\left|\frac{c_{\tau(i)}-c_{\tau(i+1)}}{2}\right|$, we have
$[{\tca^\prime}^r_{r_0}(I[\bfh(r_0)])]=[I[\bfh(r)]]$.
\end{lem}

We denote by $\hca(r)$ the set of isomorphism classes of irreducible
Harish-Chandra modules $Z$ such that $J(\bfh(r))$ annihilates $Z$ and $\Dim(Z)=\frac{1}{2}\dim\hol(\bfc)$. 
If we apply [Vogan 1990] Corollary 7.14 to our setting, we have:

\begin{lem}

We fix $r\in\nat$ such that
 $|r-h_{\tau(i)}+h_{\tau(i+1)}|\leqslant \left|\frac{c_{\tau(i)}-c_{\tau(i+1)}}{2}\right|$.
(We easily see such an $r$ exists.)

(1) \,\,\, $\tca^r_0$ gives a bijection of $\hca(0)$ onto $\hca(r)$.

(2) \,\,\, For sufficiently large $r_0\in\nat$, ${\tca^\prime}^r_{r_0}$
 gives a bijection of $\hca(r_0)$ onto $\hca(r)$.
\end{lem}

From Lemma 2.9.2, Lemma 2.9.3, and the exactness of the translation
functors,
we see that the number of irreducible constituent of the maximal
Gelfand-Kirillov dimension in $I[\bfh]$ equals the number of those in
$I[\bfh(r_0)]$.
From the assumption of the induction, ($\Diamond$) holds for
$I[\bfh(r_0)]$.
Therefore, we have ($\Diamond$) for $I[\bfh]$. \,\,\,\,\,Q.E.D.

\subsection{The socle of a degenerate principal series representation
  with a finite-dimensional quotient}

In order to obtain the main result, we prepare the following lemma in
the general setting.

\begin{lem}
Let $G$ be a real reductive group and let $K$ be a maximal compact
 subgroup of $G$.
We denote by $\gggg$ the complexified Lie algebra of $G$.
For a Harish-Chandra $(\gggg,K)$-module $M$, we denote by $M^h$ the
 Hermitian dual of $M$. (cf.\ [Vogan 1984])

Let $I$ be a Harish-Chandra $(\gggg,K)$-module and let $\Psi : I^h\rightarrow I$ be a $(\gggg,K)$-homomorphism.

We assume the following three conditions (a1)-(a3).

(a1) \,\,\, $\Dim({\rm Kernel}(\Psi))<\Dim(I^h)$.

(a2) \,\,\, Let $V$ be an arbitrary irreducible constituent of $I$ such
 that $\Dim(V)=\Dim(I)$.  Then, $V^h\cong V$ and the multiplicity of $V$
 in $I$ is one.

(a3) \,\,\, Let $W$ be any submodule of $I$, then $\Dim(W)=\Dim(I)$.

Then, we have

(1) \,\,\, ${\rm Socle}(I)=\Psi(I^h)$.

(2) \,\,\, $\Dim(I/{\rm Socle}(I))<\Dim(I)$.

Here, ${\rm Socle}(I)$ is the maximal semisimple submodule of $I$.

In particular, any irreducible constituent $V$ of $I$ such that
 $\Dim(V)=\Dim(I)$ is a submodule of $I$.
\end{lem}

\proof

Let $V$ be any irreducible submodule of $I$. 
Then $V^h\cong V$ has to be realized as a quotient of $I^h$.
From (a3), we have $\Dim(V)=\Dim(I)$.
From (a1), we see $\Psi(I^h)$ contains an irreducible constituent
isomorphic to $V$. From (a2), we have $V\subseteq\Psi(I^h)$.
Hence ${\rm Socle}(I)\subseteq\Psi(I^h)$.

In order to show (1), we have only to show that $\Psi(I^h)\cong I^h/{\rm
Kernel}(\Psi)$ is semisimple.
Let $W$ be any irreducible quotient of $I^h/{\rm
Kernel}(\Psi)$.
Then, $W$ is also a quotient of $I^h$.
So, $W^h$ is realized as a submodule of $I$.
Hence, $\Dim(W^h)=\Dim(W)=\Dim(I)$ from (a3).
So, $W\cong W^h\subseteq {\rm Socle}(I)\subseteq\Psi(I^h)$ 
Hence, $I^h/{\rm Kernel}(\Psi)$ contains a submodule isomorphic to $W$.
(a2) implies the multiplicity of $W$ is one.
Since any irreducible quotient of $I^h/{\rm Kernel}(\Psi)$ is a
submodule, $\Psi(I^h)\cong I^h/{\rm
Kernel}(\Psi)$ is semisimple. 

(2) follows from (a1) and (1).  \,\,\,\,$\blacksquare$
 
Now, we state the main result.

\begin{thm}
Assume $\bpu=(u_1,...,u_s), \bpv=(v_1,...,v_s)\in\itg^s$ and $h\in\itg$
 satisfy $u_1\geqslant\cdots\geqslant u_s\geqslant h\geqslant v_s\geqslant\cdots\geqslant v_1$.
Put
 $\bpc(m,n,\kappa)=(k_1,...,k_s,m+n-2k^\ast_s,k_s,...,k_1)\in\nat^{2s+1}$.
We denote by  $\delta^\kappa_{m,n}$ an element of $\cpx^s$ whose the
$i$-th entry is $\frac{m+n-k^\ast_{i-1}}{2}$, where $k^\ast_i=k_1+\cdots+k_i$. 

(1) \,\,\,\, We have

\[
 {\rm Socle}({}^nI^\kappa_{m,n}[\bpu+\delta^\kappa_{m,n};h;
 \bpv-\delta^\kappa_{m,n}])
\cong\bigoplus_{(\bpm,\bpn)\in\dbO(\bpc(m,n;\kappa))}\ana_{(\bpm,\bpn)}[u_1,...,u_s,h,v_s,...,v_1].
\]

(2) \,\,\,\, There exists an embedding of a generalized Verma module:
\[
 {}^nM_{\ppp_\kappa(m,n)}(-\xi(\bfu+\delta^\kappa_{m,n},h,\bfv-\delta^\kappa_{m,n}))\hookrightarrow {}^nM_{\ppp_\kappa(m,n)}(-\xi(\bfv-\delta^\kappa_{m,n},h,\bfu+\delta^\kappa_{m,n})).
\]
This embedding induces an intertwining operator 
\[
\varphi : {}^nI^\kappa_{m,n}[\bpv-\delta^\kappa_{m,n};h;
 \bpu+\delta^\kappa_{m,n}])\rightarrow {}^nI^\kappa_{m,n}[\bpu+\delta^\kappa_{m,n};h;
 \bpv-\delta^\kappa_{m,n}]).
\]
Moreover, we have ${\rm Image}(\varphi)={\rm Socle}({}^nI^\kappa_{m,n}[\bpu+\delta^\kappa_{m,n};h;
 \bpv-\delta^\kappa_{m,n}]).$

\end{thm}

\proof

From the standard argument of translation principle, we may assume that 
$u_1=u_2=\cdots=u_s=h=v_s=\cdots=v_1=0$.
For simplicity, we put $I={}^nI^\kappa_{m,n}[\delta^\kappa_{m,n};0;
 -\delta^\kappa_{m,n}])$ and $\ppp=\ppp_\kappa(m,n)$.
Then, $I^h={}^nI^\kappa_{m,n}[-\delta^\kappa_{m,n};0;
 +\delta^\kappa_{m,n}])$.
Moreover, $I$ (resp.\ $I^h$) has a perfect pairing with a generalized
Verma module ${}^uM_\ppp(-2)$ (resp.\ ${}^uM_\ppp(0)$). (See 1.2.)
In this case $\ppp$ satisfies the equivalent conditions in Proposition
1.2.1.
Hence, Proposition 1.2.7 implies that ${}^uM_\ppp(-2)\hookrightarrow
{}^uM_\ppp(0)$.
This induces an intertwining operator $\Psi : I^h\rightarrow I$ (for
instance, see [Collingwood-Shelton 1990] 2).
Since ${\rm Kernel}(\Psi)$ has a perfect pairing with
${}^uM_\ppp(0)/{}^uM_\ppp(-2)$ and
$\Dim({}^uM_\ppp(0)/{}^uM_\ppp(-2))<\Dim({}^uM_\ppp(-2))=\Dim(I)$, we
have $\Dim({\rm Kernel}(\Psi))<\Dim(I)$. 
Since ${}^uM_\ppp(-2)$ is irreducible, any submodule of $I$ also has a
parfect pairing with ${}^uM_\ppp(-2)$.
This implies that any submodule of $I$ has the same Gelfand-Kirillov
dimension as $I$.
Together with Lemma 2.10.1, we can apply Theorem 2.9.1 to our setting and
obtain the result.
\,\,\,\, Q.E.D.

The Matsuki duality ([Matsuki 1988]) tells us $K_\cpx(m,n)$-conjugacy class
of $\theta$-stable parabolic subalgberas $\{\qqq(\bpm,\bpn)\mid
\bpm+\bpn=\bpc(m,n,\kappa)\}$ parameterizes the open $\U(m,n)$-orbits of
the complexified generalized flag variety $X(m,n,\kappa)=\GL(m+n,\cpx)/P_\kappa(m,n)_\cpx$.
We denote $\lca$ the canonical line bundle on $X(m,n,\kappa)$.
For an open  $\U(m,n)$-orbit $\hol$ in $X(m,n,\kappa)$, we put 
\[
 \ana_\hol=H^S(\hol,\lca)_{\mbox{$K(m,n)$-finite}}\cong \ana_{(\bpm,\bpn)}[0,..,0].
\]
Here, $S=\dim(\vvv(\bpm,\bpn)\cap\kkk(m,n))$ and $\qqq(\bpm,\bpn)$ is
the $\theta$-stable parabolic subalgebra corresponding to $\hol$.
We call an open  $\U(m,n)$-orbit $\hol$ in $X(m,n,\kappa)$ good, if
$\Dim(\ana_\hol)=\dim X(m,n,\kappa)$.
Then we may rewrite Theorem 2.10.2 as follows.

\begin{cor}
\[
{\rm Socle}({}^u\Ind_{P_\kappa(m,n)}^{\U(m,n)}(\lca))\cong
 \bigoplus_{\mbox{$\hol$ : a good open $\U(m,n)$-orbit in
 $X(m,n,\kappa)$}}\ana_\hol.
\]
\end{cor}


\setcounter{section}{3}
\section*{\S\,\, 3.\,\,\,\, $\GL(n,\rel)$  and $\GL(n,\qtr)$}
\setcounter{subsection}{0}

\subsection{{Parabolic subalgebras and generalized flag varieties}}

Let $n$ be a positive integer greater than $1$.
Let $\dbK$ be one of $\rel$, $\cpx$, or $\qtr$.
We consider $G=\GL(n,\dbK)$.
We denote by $\sH$ the Lie group of $n\times n$ diagonal
non-singular matrices in $G$.
Fix  a maximal compact subgroup $K$ of $G$.
We denote by $\theta$ the corresponding Cartan involution.
We may choose the maximal compact subgroup $K$ appropriately so that
$\sH(n)$ is $\theta$-stable.
We denote by $\gggg$ (resp.\ $\kkk$) the complexified Lie algebra of $G$
(resp.\ $K$).
We denote by  $\shh$ the  complexified Lie algebra of $\sH$.
Then, $\shh$ is a  maximally split Cartan subalgebra
of $\gggg$. 

Let $\bfc=(c_1,...,c_s)$ be a sequence of positive integers such that
$c_1+\cdots+c_s=n$.
We denote by $P(\bfc)$ the block-upper-triangular parabolic subgroup of $G$
with blocks of sizes $c_1,...,c_s$ along the diagonal.
We denote by $P_\cpx(\bfc)$ (resp.\ $G_\cpx$) the complexification of
$P(\bfc)$ (resp.\ $G$).
We denote by $\ppp(\bfc)$ the complexified Lie algebra of $P(\bfc)$.
We denote by $\mmm(\bfc)$ the Levi subalgebra of $\ppp(\bpc)$
such that $\shh(n)\subseteq\mmm(\bfc)$.
and denote by $M_\cpx(\bfc)$ the corresponding complexified Levi subgroup.
We denote by $\nnn(\bfc)$ the nilradical of $\ppp(\bpc)$.

We denote by $w_0$ (resp.\ $w_{\bfc}$) the longest element of the Weyl
group for $(\gl(n,\cpx),\shh)$ (resp.\ $(\mmm(\bfc),\shh)$).

From [Matsuki 1979], we easily see the full flag variety of
$G_\cpx$ has a unique open $G$-orbit.
Hence, an arbitrary generalized flag variety $X(\bfc)=G_\cpx/P_\cpx(\bfc)$
has a  unique open $G$-orbit (say $\hol_{\bfc}$).
We regard as $X(\bfc)$ the set of parabolic subalgebras of $\gggg$
which are $\Ad(G_\cpx)$-conjugate to $\ppp(\bfc)$.
The following result is more or less known.

\begin{prop}
Let $\bfc=(c_1,...,c_s)$ be a sequence of positive integers such that
$c_1+\cdots+c_s=n$.
Then, the following conditions are equivalent to each other.

(1) \,\,\, $\hol_{\bfc}$ contains a $\theta$-stable parabolic
subalgebra of $\gggg$ .

(2) \,\,\, $\ppp(\bfc)$ and the parabolic subalgebra opposite to $\ppp(\bfc)$ are $G$-conjugate.

(3) \,\,\, For any integer $1\leqslant i \leqslant \frac{s}{2}$, we have
 $c_i=c_{s-i+1}$.

(4) \,\,\, $w_0 w_{\bfc}=w_{\bfc} w_0$. (Namely, $w_0w_{\bfc}$ is an involution.)
\end{prop}

We denote by $\omega_{\bfc}$ the canonical line bundle on $X(\bfc)$.
Namely, $\omega_{\bfc}$ is a homogeneous line bundle associate with the
character of $M_\cpx(\bfc)$ : 
\[
 M_\cpx(\bfc)\ni g\rightsquigarrow \det\left(\Ad(g)|_{\nnn(\bfc)}\right)\in\cpx^\times.
\]
 
\subsection{The socle of degenerate principal series representation with
respect to canonical bundle}

Assume that $\bfc$ satisfies the equivalent conditions in Proposition
3.1.1 and choose a $\theta$-stable parabolic subgroup $\qqq(\bfc)$ in
$\hol_{\bfc}$ such that $\qqq(\bfc)$ has a Levi part $\llll(\bfc)$ which
is a comlexified Lie algebra of its normalizer $L(\bfc)$ in $G$. 
We denote by $\uuu(\bfc)$ the nilradical of $\qqq(\bfc)$

Put
$\ana_{\hol_{\bfc}}=H^S(\hol_{\bfc},\omega_{\bfc})_{\mbox{$K$-finite}}\cong
\ana_{\qqq(\bfc)}(0).$
Here, $S=\dim(\kkk\cap\uuu(\bfc))$.

The case of $\dbK=\cpx$ is discussed in Corollary 1.2.8.
For the remaining cases, we have:

\begin{thm} \,\,\,\,
\mbox{}

(1) \,\,\,\,  If $\dbK=\qtr$, then $\ana_{\hol_{\bfc}}$ is a unique
 irreducible submodule of ${}^u\Ind_{P(\bfc)}^G(\omega_{\bfc})$.

(2) \,\,\,\, If $\dbK=\rel$, then there is a one-dimensional
 representation $\chi$ which  is trivial on the identical component of
 $P(\bfc)$ and $\ana_{\hol_{\bfc}}$ is a unique
 irreducible submodule of
 ${}^u\Ind_{P(\bfc)}^G(\omega_{\bfc}\otimes\chi)$.

(3) \,\,\,\, In any case, $\ana_{\hol_{\bfc}}$ is the image of an
 intertwining operator induced from an embedding of a generalized Verma
 module.

\end{thm}

\proof

We consider the case of  $\dbK=\qtr$.
Write $\qrt=\cpx+j\cpx$.
We put 
$K=\Spp(n)=\{g\in \GL(k,\qrt)| {}^t\bar{g}g=I_k \}$.
Then we regard $\gl(k,\cpx)$ as a real Lie subalgebra of $\gl(k,\qrt)$.
For $\ell\in\itg$ and $t\in\sqrt{-1}\rel$, we define a one-dimensional unitary representation
$\xi_{\ell}$ of $\GL(k,\cpx)$ as follows.
\[
 \xi_{\ell}(g)=\left(\frac{\det(g)}{|\det(g)|}\right)^\ell.
\]
Let $\qqq(k)$ be a $\theta$-stable  parabolic subalgebra with a Levi decomposition $\qqq(k)=\gl(k,\cpx)+\uuu(k)$.
We choose the nilradical $\uuu(k)$ so that $\xi_{\ell}$ is in the good range with respect
to $\qqq(k)$ for sufficient large $\ell$.
Derived functor modules with respect to $\qqq(k)$ is called quarternionic Speh representations.
\begin{equation}
 A_k(\ell)=({}^u\rca_{\qqq(k),\Oo(k)}^{\gl(k,\qrt)\otimes_\rel\cpx,
\Spp(k)})^{k(k+1)}(\xi_{\ell+2k})   \,\,\,\,\,  \,\, (\ell \in\itg).\tag{$\ast$}
\end{equation}
For $\ell\in\itg$, $A_k(\ell)$ is derived functor module in the good (resp.\ weakly fair) range
in the sense of [Vogan 1988] if and only if $\ell\geqslant 0$ (resp.\
$\ell\geqslant -k$).

We define $\bfd=(d_1,...,d_{[\frac{s+1}{2}]} )$ by
$d_{i}=2c_i=2c_{s-i+1}$ \,\,\,\, $(1\leqslant
i\leqslant\frac{s}{2})$.
Here, $[\frac{s+1}{2}]=\max\{j\in\itg\mid j\leqslant\frac{s+1}{2}\}$.
If $s$ is odd, we put $d_{\frac{s+1}{2}}=c_{\frac{s+1}{2}}$.
We put $d^\ast_i=d_1+\cdots+d_i$.
The Levi part $M(\bfd)$ of $P(\bfd)$ can be identified with
$\GL(d_1,\qtr)\times\cdots\times\GL(d_{[\frac{s+1}{2}]},\qtr)$. 

For $g\in GL(k,\qtr)$, we denote by $\det(g)\in\rel^\times$ the noncommutative
determinant ([Artin 1957] p151).
For $t\in\rel$, we denote by $\chi_t$ a character of $\GL(k,\qtr)$
defined by $\chi_t(g)=\det(g)^t$.

More or less [Vogan 1986] tells us (also see [Matumoto 2002] 2.4 )
\begin{equation*}
 \ana_{\hol_{\bfc}}\cong
 {}^n\Ind_{P(\bfd)}^G\left(\bigotimes_{i=1}^{[\frac{s}{2}]}
 A_{d_i}(2n-2d^\ast_i)
 \boxtimes id\right).  \tag{$\bigstar$}
\end{equation*}

Here, $id$ means the identity representation of
$\GL(n-d_{[\frac{s}{2}]}^\ast,\qtr)$.
(We regard $\GL(0,\qtr)$ as the trivial group.)
``$\bigotimes$''in the above formula ($\bigstar$) means an external tensor product.
We quote:

\begin{lem} \,\,\,\, ([Sahi 1995], [Zhang 1995])
\begin{equation*}
 A_{2k}(\ell)\sim {}^n\Ind_{P((k,k))}^{\GL(2k,\qtr)}(\chi_{\ell+k}\boxtimes\chi_{-\ell-k}). 
\end{equation*}

\end{lem}

We define $\bfc^\prime=(c_1^\prime,...,c_s^\prime)$ by
$c_{2i-1}^\prime=c_{2i}^\prime=c_i=c_{s-i+1}$ \,\,\,\, $(1\leqslant
i\leqslant\frac{s}{2})$.
If $s$ is odd, we put $c_s=c_{\frac{s+1}{2}}$.
Combining the above ($\bigstar$) and Lemma 3.2.2, we can easily see
\[
 \ana_{\hol_{\bfc}}\sim
 {}^n\Ind_{P(\bfc^\prime)}^G\left(\bigotimes_{i=1}^{[\frac{s}{2}]}\left(\chi_{2n-2d^\ast_i+\frac{d_i}{2}}\boxtimes\chi_{-2n+2d^\ast_i-\frac{d_i}{2}}\right)\boxtimes
 id\right).
\]
Since $M(\bfc)$ is $G$-conjugate to $M(\bfc^\prime)$, Harish-Chandra's
character theory implies the following character identity.
\[
 \left[ {}^n\Ind_{P(\bfc^\prime)}^G\left(\bigotimes_{i=1}^{[\frac{s}{2}]}\left(\chi_{2n-2d^\ast_i+\frac{d_i}{2}}\boxtimes\chi_{-2n+d^\ast_i-\frac{d_i}{2}}\right)\boxtimes
 id\right)\right]=[{}^u\Ind_{P(\bfc)}^G(\omega_{\bfc})].
\] 
Hence,
\[
 \ana_{\hol_{\bfc}}\sim {}^u\Ind_{P(\bfc)}^G(\omega_{\bfc}).
\]
We can apply an argument similar to the proof of Theorem 2.10.2 and we
have the desired result for the case of $G=\GL(n,\qtr)$.

For the case of $G=\GL(n,\rel)$, we can prove the result  similarly.
We may apply [Speh 1983] instead of ($\bigstar$).
The corresponding result of result of Lemma 3.2.2 is also obtained in
[Sahi-Stein 1990], [Sahi 1995], and  [Zhang 1995]. (Also see [Howe-Lee 1999])
\,\,\,\, $\blacksquare$


\bigskip

\subsection*{References}

\mbox{}\hspace{5mm}

{\sf [Artin 1957]} E.\ Artin, Geometric algebra,
Interscience Publishers, Inc., New York-London, 1957. 

{\sf [Barbasch-Vogan 1980]} D.\ Barbasch and D.\ A.\ Vogan Jr.\ , The local structure of characters, {\it J.\ Funct.\ Anal.\ } {\bf 37} (1980), 27-55.

{\sf [Barbasch-Vogan 1983]} D.\ Barbasch and D.\ A.\ Vogan Jr., Weyl group representations and nilpotent orbits, in: P.\ C.\ Trombi, editor, ``Representation Theory of Reductive Groups,'' Progress in Mathematics Vol.\ 40, 21-33, Birkh\"{a}user, 
Boston-Basel-Stuttgart, 1983.

{\sf [Bernstein-Gelfand 1980]} J.\ Bernstein and S.\ I.\ Gelfand, Tensor product of finite and infinite dimensional representations of semisimple Lie algebras, {\it Compos.\ Math.\ } {\bf 41} (1980), 245-285.

{\sf [Borho-Jantzen 1977]} W.\ Borho and J.\ C.\ Jantzen, \"{U}ber primitive Ideale in der Einh\"{u}llenden einer halbeinfachen Lie-Algebra, {\it Invent.\ Math.\ } {\bf 39} (1977), 1-53.

{\sf [Collingwood-Shelton 1990]} D.\ H.\ Collingwood and B.\ Shelton,
A duality theorem for extensions of induced highest weight modules,
{\it Pacific J.\ Math.\ } {\bf 146} (1990), 227--237.

{\sf [Conze-Berline, Duflo 1977]} N.\ Conze-Berline and M.\ Duflo, Sur
les repr\'{e}sentations induites des groupes semi-simples complexes,
{\it Comp.\ Math.\ } {\bf 34} (1977), 307-336.

{\sf [Duflo 1977]} M.\ Duflo, Sur la classifications des id\'{e}aux primitifs dans l'alg\`{e}bre de Lie semi-simple, {\it Ann.\ of Math.\ } {\bf 105} (1977), 107-120.

{\sf [Enright 1979]} T.\ J.\ Enright, Relative Lie algebra cohomology and
unitary representations of complex Lie groups, {\it Duke Math.\ J.\ }{\bf 46} (1979), 513-525.

{\sf [Gindikin 1993]} S.\ Gindikin, 
Holomorphic language for $\overline\partial$-cohomology and
representations of real semisimple Lie groups, in : 
The Penrose transform and analytic cohomology in representation theory (South Hadley, MA, 1992), 103--115,
Contemp. Math., 154,
Amer. Math. Soc., Providence, RI, 1993.

{\sf [Howe-Lee 1999]} R.\ Howe, and S.\ T.\ Lee, 
Degenerate principal series representations of ${\rm GL}\sb n(\bold C)$ and ${\rm GL}\sb n(\bold R)$, 
{\it J.\ Funct.\ Anal.\ } {\bf 166} (1999), 244--309.

{\sf [Jantzen 1983]} J.\ C.\ Jantzen, ``Einh\"{u}llende Algebren halbeinfacher Lie-Algebren'' Ergebnisse der Mathematik und ihrer Grenzgebiete 3, Springer-Verlag, Berlin, Heidelberg, New York, Tokyo, 1983.

{\sf [Joseph 1979]} A.\ Joseph, Dixmier's problem for Verma and principal series submodules, {\it J.\ London Math.\ Soc.\ } (2) {\bf 20} (1979), 193-204.

{\sf [Joseph 1983]} A.\ Joseph,
On the classification of primitive ideals in the enveloping algebra of a semisimple Lie algebra.
in: Lie group representations, I (College Park, Md., 1982/1983), 30--76,
Lecture Notes in Math., 1024,
Springer, Berlin, 1983.

{\sf [Kashiwara-Kawai-Kimura 1986]}
M.\ Kashiwara, T.\ Kawai, and T.\ Kimura, 
Foundations of algebraic analysis.
Translated from the Japanese by Goro Kato. Princeton Mathematical Series, 37.
Princeton University Press, Princeton, NJ, 1986.

{\sf [Kashiwara-Laurent 1983]} M.\ Kashiwara and Y.\ Laurent,
Th\'{e}or\`{e}mes d'annulation et deuxi\`{e}me microlocalisation,
Universit\`{e} de Paris-Sud, 1983.

{\sf [Kashiwara-Schapira 1990]} M.\ Kashiwara and P.\ Schapira, 
Sheaves on manifolds,
Grundlehren der Mathematischen Wissenschaften, 292,
Springer-Verlag, Berlin, 1990.

{\sf [Kostant 1961]} B.\ Kostant, 
Lie algebra cohomology and the generalized Borel-Weil theorem,
{\it Ann.\ of Math.\ } {\bf 74} (1961) 329--387.

{\sf [Matsuki 1979]} T.\ Matsuki,
The orbits of affine symmetric spaces under the action of minimal parabolic subgroups,
{\it J.\ Math.\ Soc.\ Japan} {\bf 31} (1979), 331--357.

{\sf [Matsuki 1988]} T.\ Matsuki, Closure relations for orbits on
affine symmetric spaces under the action of parabolic subgroups,
Intersections of associated orbit {\it Hiroshima Math.\ J.\ } {\bf 18} (1988), 59-67.

{\sf [Matumoto 1988]} H.\ Matumoto, 
Cohomological Hardy space for ${\rm SU}(2,2)$,
in : Representations of Lie groups, Kyoto, Hiroshima, 1986, 469--497,
Adv. Stud. Pure Math., 14,
Academic Press, Boston, MA, 1988.
 
{\sf [Matumoto 1993]} H.\ Matumoto, On the existence of homomorphisms
between scalar generalized Verma modules, in: Contemporary
Mathematics, 145, 259-274, Amer. Math. Soc., Providence, RI,
1993. 

{\sf [Matumoto 1996]} H.\ Matumoto, On the representations of $U(m,n)$
unitarily induced from derived functor modules, {\it Compositio Math.\ }
{\bf 100} (1996), 1-39.

{\sf [Matumoto 2002]} H.\ Matumoto, On the representation of Sp(p,q) ans SO*(2n) 
unitarily induced from derived functor modules, to appear in {\it
Compositio Math.\ }. arXive : math.RT/0203107

{\sf [Sahi 1992]} S.\ Sahi,
The Capelli identity and unitary representations,
{\it Compositio Math.\ } {\bf 81} (1992), 247--260.

{\sf [Sahi 1995]} S.\ Sahi,
 Jordan algebras and degenerate principal series.
{\it J.\ Reine Angew.\ Math.\ } {\bf 462} (1995), 1--18.

{\sf [Sahi-Stein 1990]} S.\ Sahi and E.\ M.\ Stein, 
Analysis in matrix space and Speh's representations,
{\it Invent.\ Math.\ } {\bf 101} (1990), 379--393.

{\sf [Schmid-Vilonen 2000]} W.\ Schmid and K.\ Vilonen, 
Characteristic cycles and wave front cycles of representations of reductive Lie groups,
{\it Ann.\ of Math.\ } {\bf 151} (2000), 1071--1118.

{\sf [Speh 1983]} B.\ Speh, Unitary representations of $\GL(n,\rel)$ with
non-trivial $(\gggg, K)$-cohomology, {\it Invent Math.\ }{\bf 71} (1983), 443-465.

{\sf [Vogan 1982]} D.\ A.\ Vogan Jr., ``Representations of Real reductive Lie Groups'', {\it Progress in Mathematics,} Birkh\"{a}user, 1982.

{\sf [Vogan 1984]} D.\ A.\ Vogan Jr., Unitarizability of certain series of representations, {\it Ann.\ of Math.\ } {\bf 120} (1984), 141-187.

{\sf [Vogan 1986]} D.\ A.\ Vogan Jr., The orbit method and primitive
ideals for semisimple Lie algebras, in :
Lie algebras and related topics (Windsor, Ont., 1984), p281-316,
CMS Conf. Proc., 5,
Amer. Math. Soc., Providence, RI, 1986. 

{\sf [Vogan 1986]} D.\ A.\ Vogan Jr., The unitary dual of $GL(n)$ over an archimedean field, {\it Invent.\ Math.\ } {\bf 83} (1986), 449-505.

{\sf [Vogan 1988]} D.\ A.\ Vogan Jr., Irreducibilities of discrete series representations for semisimple symmetric spaces, {\it Adv.\  Stud.\ in Pure Math.\ } vol.\ 14, Kinokuniya Book Store, 1988, 381-417.

{\sf [Vogan 1990]} D.\ A.\ Vogan Jr.,
Dixmier algebras, sheets, and representation theory, in:
Operator algebras, unitary representations, enveloping algebras, and invariant theory (Paris, 1989), 333--395, Progr. Math., 92,
Birkh\"{a}user Boston, Boston, MA, 1990. 

{\sf [Vogan 1997]}\,\, D.\ A.\ Vogan Jr.,  Cohomology and group representations,
in ``Representation theory and automorphic forms (Edinburgh, 1996)'', 
Proc. Sympos. Pure Math., 61. 219-243, Amer. Math. Soc., Providence, RI,
1997.

{\sf [Vogan-Zuckerman 1984]} D.\ A.\ Vogan Jr.\ and G.\ Zuckerman,
Unitary representations with non-zero cohomology, {\it Compositio Math.\
}{\bf 53} (1984), 51-90.

{\sf [Wong 1993]} H.\ Wong, 
Dolbeault cohomologies and Zuckerman modules, in :
The Penrose transform and analytic cohomology in representation theory (South Hadley, MA, 1992), 217--223,
Contemp.\ Math.\, 154,
Amer.\ Math.\ Soc.\, Providence, RI, 1993. 

{\sf [Zhang 1995] }  G.\ K.\ Zhang, 
Jordan algebras and generalized principal series representations,
{\it Math.\ Ann.\ } {\bf 302} (1995), 773--786.

\end{document}